\documentclass[aos,MSNbibl,nameyear,dvips]{arximspdf}
\usepackage{graphicx}
\doi{10.1214/14-AOS1280} % kopijuoti is laisko
\volume{43}
\issue{1}
\pubyear{2015}
\firstpage{276}
\lastpage{298}
\docsubty{FLA}

\makeatletter
\newcommand{\rrvert}{\vert}
\newcommand{\llvert}{\vert}
\newcommand{\eqref}[1]{(\ref{#1})}

\newtheorem{lem}{Lemma}[section]
\newtheorem{theo}{Theorem}[section]

\newproclaim{rem}{Remark}[section]
\newproclaim{defi}{Definition}[section]
\newproclaim{ex}{Example}[section]
\newproclaim{cond}{Condition}[section]

\newcommand{\R}{\mathbb{R}}

\makeatother

\begin{document}
\begin{frontmatter}

\title{On the block maxima method in extreme value theory: PWM estimators\thanksref{T1}}
\runtitle{On the block maxima method}
\thankstext{T1}{Supported in part by FCT Project PTDC /MAT /112770 /2009;
EXPL/MAT-STA/0622/2013 and PEst-OE/MAT/UI0006/2014.}

\begin{aug}
\author[A]{\fnms{Ana} \snm{Ferreira}\corref{}\ead[label=e1]{anafh@isa.utl.pt}}
\and
\author[B]{\fnms{Laurens} \snm{de Haan}\ead[label=e2]{ldehaan@ese.eur.nl}}
\address[A]{ISA\\
Department of Economics\\
University of Lisbon\\
Tapada da Ajuda 1349-017 Lisboa\\
Portugal\\
and\\
CEAUL\\
FCUL\\
Bloco C6 - Piso 4 Campo Grande\\
749-016 Lisboa\\
Portugal\\
\printead{e1}}
\address[B]{Department of Mathematics\\
Erasmus University Rotterdam\\
P.O. Box 1738\\
3000 DR Rotterdam\\
The Netherlands\\
and\\
CEAUL\\
FCUL\\
Bloco C6 - Piso 4 Campo Grande\\
749-016 Lisboa\\
Portugal\\
\printead{e2}}
\runauthor{A. Ferreira and L. de Haan}
% \and
\affiliation{University of Lisbon and Erasmus Univ Rotterdam}
%%\address[A]{\\\printead{e1}}
\end{aug}

% HISTORY:
\received{\smonth{4} \syear{2014}}
\revised{\smonth{8} \syear{2014}}

% ABSTRACT

\begin{abstract}
In extreme value theory, there are two fundamental approaches, both
widely used: the block maxima (BM) method and the peaks-over-threshold
(POT) method. Whereas much theoretical research has gone into the POT
method, the BM method has not been studied thoroughly. The present
paper aims at providing conditions under which the BM method can be
justified. We also provide a theoretical comparative study of the
methods, which is in general consistent with the vast literature on
comparing the methods all based on simulated data and fully parametric
models. The results indicate that the BM method is a rather efficient
method under usual practical conditions.

In this paper, we restrict attention to the i.i.d. case and focus on
the probability weighted moment (PWM) estimators of
Hosking, Wallis and Wood [\textit{Technometrics} (1985) \textbf{27} 251--261].
\end{abstract}

% KEYWORDS

\end{frontmatter}

%s1 #&#
\section{Introduction}\label{Introd_sect}
The \textit{block maxima} (BM) approach in extreme value theory (EVT),
consists of dividing the observation period into nonoverlapping
periods of equal size and restricts attention to the maximum
observation in each period [see, e.g., \citet{Gumbel58}]. The new
observations thus created follow---under domain of attraction
conditions, cf. \eqref{mda_cond} below---approximately an extreme
value distribution, $G_\gamma$ for some real $\gamma$. Parametric
statistical methods for the extreme value distributions are then
applied to those observations. %The procedure can be justified using (a
%strengthening of) the domain of attraction conditions of EVT.

In the \textit{peaks-over-threshold} (POT) approach in EVT, one selects
those of the initial observations that exceed a certain high threshold.
The probability distribution of those selected observations is
approximately a generalized Pareto distribution \citet{Pickands75}.
%%Parametric statistical methods for GPD are then applied to those
%observations. Again, a bias may appear since GPD is not the exact
%distribution of those selected observations.

%The BM method is the older one. The POT method has been developed by
%statistical tools.

In the case of the POT method, exact conditions under which the
statistical method is justified can be described by a second-order term
[see, e.g., \citet{Drees98} and \citet{HaanFerreira}, Section~2.3].
In the case of block maxima, usually it is taken for granted that the
maxima follow very well an extreme value distribution. In this paper,
we take this misspecification into account by quantifying it in terms
of a second-order expansion; cf. Condition \ref{2ndordcond} below.
Since $G_\gamma$ is not the exact distribution for those observations,
a bias may appear.
%, and it is one of the aims of this paper to formulate exact
%conditions under which the statistical method can be justified. %For
%instance, since some of the block maxima may actually not be very
%high, one expects somewhat more strict conditions in this case then in
%the POT case. However, as it turns out, the conditions are similar.

The POT method picks up all ``relevant'' high observations. The BM method
on the one hand misses some of these high observations and, on the
other hand, might retain some lower observations. Hence the POT seems
to make better use of the available information.
%Indeed this is in agreement with several studies based on simulated
%data and fully parametric models: e.g.
%the general opinion among statisticians which have been seen
%established over the recent years.
%E.g. in Katz, Parlange and Naveau (2002) one can read `Its rationale
%is that if additional information about the extreme upper tail were
%used besides the annual maxima (i.e., other relatively high values in
%the sample), then more accurate estimates of the parameters and
%quantiles of extreme value distributions would be obtained';

There are practical reasons for using the BM method:
\begin{itemize}
\item The only available information may be block maxima, for example,
yearly maxima with long historical records or long range simulated data
sets \citet{Kharin07}.
\item The BM method may be preferable when the observations are not
exactly independent and identically distributed (i.i.d.). For example,
there may be a seasonal periodicity in case of yearly maxima or, there
may be short range dependence that plays a role within blocks but not
between blocks; cf. for example, Katz, Parlange and Naveau (\citeyear{KatzParlangeNaveau02}) and
\citet{MadsenRasmussenRosbjerg97} for further discussion.
\item The BM method may be easier to apply since the block periods
appear naturally in many situations [\citet{NaveauGuillouCooleyDiebolt09},
van~den Brink, K{\"{o}}nnen and Opsteegh (\citeyear{BrinkKonnenOpsteegh05}),
\citet{Valk93}]. On the other hand, the POT method allows
for greater flexibility in many cases since it might be difficult to
change the block size in practice. %Moreover the problem of choosing a
%high threshold in the POT method (which is a difficult one) does not
%play a role ().
\end{itemize}

When working with BM, there are two sets of estimators that are widely
used: the maximum likelihood (ML) estimators [e.g., \citet{PrescottWalden1980}] and the probability weighted moment (PWM) estimators
\citet{HoskingWallisWood85}. Recently, \citet{Dombry} has proved
consistency of the former. The present paper concentrates on the
latter. Our work has given rise to the paper \citet{BucherSegers14}
on the multivariate case.

The PWM estimators under the $G_\gamma$ model are very popular, for
example, in applications to hydrologic and climatologic extremes,
because of their computational simplicity, good performance for small
sample sizes and robustness even for location and scale parameters
[\citet{DieboltGuillouNaveauRibereau08},
Katz, Parlange and Naveau (\citeyear{KatzParlangeNaveau02}),
\citet{Caires09,Hosking90}].

The relative merits of POT and BM have been discussed in several
papers, all based on simulated data: \citet{Cunnane73} states that
for $\gamma=0$ and ML estimators, the POT estimate for a high quantile
is better only if the number of exceedances is larger than 1.65 times
the number of blocks; \citet{Wang91} writes that POT is as
efficient as BM model for high quantiles, based on PWM estimators;
Madsen, Pearson and Rosbjerg (\citeyear{MadsenPearsonRosbjerg97})
and
Madsen, Rasmussen and Rosbjerg (\citeyear{MadsenRasmussenRosbjerg97}) write that POT is
preferable for $\gamma>0$, whereas for $\gamma<0$, BM is more
efficient, again with the number of exceedances larger than the number
of blocks; \citet{MartinsStedinger} state that the gains
(when using historical data) with the BM model are in the range of the
gains with the POT model, based on ML estimators; \citet{Caires09}
in a vast simulation study writes that with POT samples having an
average of two or more observations per block, the estimates are more
accurate than the corresponding BM estimates, and with more than 200
years of data the accuracies of the two approaches are similar and
rather good, based on several estimators including the PWM and ML estimators.

From all these studies, some even with mixed views, the following two
features seem dominant. First, POT is more efficient than BM in many
circumstances, though needing, on average, a number of exceedances
larger than the number of blocks. Secondly, POT and BM often have
comparable performances, for example, for large sample sizes.

Our theoretical comparison shows that BM is rather efficient. The
asymptotic variances of both extreme value index and quantile
estimators are always lower for BM than for POT. Moreover, the
approximate minimal mean square error is also lower for BM under usual
circumstances. The optimal number of exceedances is generally higher
than the optimal number of blocks.

The paper is organized as follows. In Section~\ref{AsympN_sect}, we
state exact conditions to justify the BM method, along with the
asymptotic normality result for the PWM estimators including high
quantile estimators. In Section~\ref{blockmaxPOT_sect}, we provide a
theoretical comparison between the two methods, BM and POT. The
analysis is based on a uniform expansion of the relevant quantile
process given in Section~\ref{Asymptnorm_sect}. This expansion also
provides a basis for analysing alternative estimators besides the PWM
estimator. Proofs are postponed to Section~\ref{proofs_sect}.

Throughout the paper, we assume that the observations are i.i.d.
In future work, we shall extend the results to the non-i.i.d. case and
to the maximum likelihood estimator.

%s2 #&#
\section{The estimators and their properties}
\label{AsympN_sect}

Let $\tilde X_1, \tilde X_2, \ldots$ be i.i.d. random variables with
distribution function $F$. Define for $m=1,2,\ldots$ and $i=1,2,\ldots,k$ the block maxima
%
%e1 #&#
\begin{equation}
\label{blockmaxima_def} X_i=\max_{(i-1)m<j\leq im} \tilde
X_j.
\end{equation}
Hence, the $m\times k$ observations are divided into $k$ blocks of size
$m$. Write $n=m\times k$, the total number of observations. We study
the model for large $k$ and $m$, hence we shall assume that $n\to\infty
$; in order to obtain meaningful limit results, we have to require that
both $m=m_n\to\infty$ and $k=k_n\to\infty$, as $n\to\infty$.

The main assumption is that $F$ is in the domain of attraction of some
extreme value distribution
\[
G_\gamma(x)=\exp \bigl( - (1+\gamma x)^{-1/\gamma} \bigr),\qquad \gamma \in
\R, 1+\gamma x>0,
\]
%
%Then $X_i$, when normalized, will follow approximately this extreme
%value distribution
that is, for appropriately chosen $a_m>0$ and $b_m$ and all $x$
%
%e2 #&#
\begin{eqnarray}
\label{mda_cond} \lim_{m\to\infty} P \biggl( \frac{X_i-b_m}{a_m}\leq x
\biggr)&=&\lim_{m\to
\infty} F^m (a_mx+b_m
)
\nonumber
\\[-8pt]
\\[-8pt]
\nonumber
&=&G_\gamma(x),\qquad i=1,2,\ldots,k.
\end{eqnarray}
This can be written as
\[
\lim_{m\to\infty} \frac{1} m \frac{1}{-\log F (a_mx+b_m )}= (1+\gamma
x)^{1/\gamma},
\]
which is equivalent to the convergence of the inverse functions:
\[
\lim_{m\to\infty}\frac{V(mx)-b_m}{a_m}=\frac{x^\gamma-1}{\gamma},\qquad x>0,
\]
with $V= (-1/\log F )^{\leftarrow}$. Hence, $b_m$ can be
chosen to be $V(m)$. This is the first-order condition. For our
analysis, we also need a second-order expansion as follows.

%co2.1 #&#
\begin{cond}[(Second-order condition)]\label{2ndordcond}
Suppose that for some positive function $a$ and some positive or
negative function $A$ with $\lim_{t\to\infty}A(t)=0$,
\[
\lim_{t\to\infty}\frac{ {(V(tx)-V(t))}/{a(t)}-{(x^\gamma-1)}/{\gamma
} }{A(t)}=\int_1^x
s^{\gamma-1}\int_1^s u^{\rho-1} \,du \,ds
=H_{\gamma,\rho}(x),
\]
for all $x>0$ [see, e.g., de Haan and Ferreira (\citeyear{HaanFerreira}),
Theorem B.3.1]. Note
that the function $|A|$ is regularly varying with index $\rho\leq0$.
\end{cond}

Let $X_{1,k},\ldots,X_{k,k}$ be the order statistics of the block
maxima $X_1,\ldots,X_k$. The statistics $\beta_0=k^{-1}\sum_{i=1}^k
X_{i,k}$ and
%
%e3 #&#
\begin{equation}
\beta_r=\frac{1} k \sum_{i=1}^k
\frac{(i-1)\cdots(i-r)}{(k-1)\cdots
(k-r)} X_{i,k},\qquad r=1,2,3,\ldots, k>r, %r>0??
\end{equation}
are unbiased estimators of $EX_1F^{rm}(X_1)$
[\citet{LandwehrMatalasWallis79}]. The PWM estimators for $\gamma$, as well as the location
$b_m$ and scale $a_m=a([m])$, are simple functionals of $\beta_0$,
$\beta_1$ and $\beta_2$. The estimator $\hat\gamma_{k,m}$ for $\gamma$
is defined as the solution of the equation
\begin{eqnarray}
\frac{3^{\hat\gamma_{k,m}}-1}{2^{\hat\gamma_{k,m}}-1}&=&\frac{3\beta
_2-\beta_0}{2\beta_1-\beta_0},\nonumber
\\
\quad\hat a_{k,m}&=&\frac{\hat\gamma_{k,m}}{2^{\hat\gamma_{k,m}}-1} \frac
{2\beta_1-\beta_0}{\Gamma (1-\hat\gamma_{k,m} )}\quad \mbox{and}\quad
%and the estimator $\hat b_{k,m}$ of $b_m$ is
\hat b_{k,m}=\beta_0+\hat
a_{k,m}\frac{1-\Gamma (1-\hat\gamma
_{k,m} )}{\hat\gamma_{k,m}},
\end{eqnarray}
where $\Gamma(x)=\int_0^\infty t^{x-1} e^{-t} \,dt$, $x>0$
[\citet{HoskingWallisWood85}].
%For $\hat\gamma_{k,m}=0$ the estimators follow by continuity.
%Note that $-\Gamma'(1)=0.577216$ is Euler's constant.
The rationale behind the estimator of $\gamma$ becomes clear when
checking the statement of Theorem~\ref{AsymptNgab} below.

%Clearly, the given estimators are quite different from the ones in the
%POT case.

%There are other variants for $M_r$, e.g. $\frac1 k \sum_{i=1}^k \left(

%Next we discuss the conditions on the distribution function $F$. The
%domain of attraction condition,

%s2.1 #&#
\subsection{Asymptotic normality}
\label{Asymptnorm_sect}

The following theorem is the basis for\break analysing estimators in the BM
approach. Let $\lceil u \rceil$ represent the smallest integer
larger than or equal to $u$.

%th2.1 #&#
\begin{theo}\label{sup_lem} Assume that $F$ is in the domain of
attraction of an extreme value distribution $G_\gamma$ and that
Condition \ref{2ndordcond} holds. Let $m=m_n\to\infty$ and $k=k_n\to
\infty$ as $n\to\infty$, in such a way that $\sqrt k A(m)\to\lambda\in\R
$. Let $0<\varepsilon<1/2$ and $ \{ X_{i,k} \}_{i=1}^k$ be the
order statistics of the block maxima $X_1,X_2,\ldots,X_k$. Then, with
$\{E_k\}_{k\geq1}$ an appropriate sequence of Brownian bridges,
\begin{eqnarray*}
&&\sqrt{k} \biggl( \frac{X_{\lceil ks\rceil,k}-b_m}{a_0(m)}- \frac{(-\log
s)^{-\gamma}-1}{\gamma} \biggr)
\\
&&\qquad=\frac{E_k(s)}{s(-\log s)^{1+\gamma}}+\sqrt k A_0(m) H_{\gamma,\rho} \biggl(
\frac{1}{-\log s} \biggr) \\
&&\qquad\quad{}+ \bigl( s^{-1/2-\varepsilon}(1-s)^{-1/2-\gamma-\rho-\varepsilon} \bigr)
o_P(1),
\end{eqnarray*}
as $n\to\infty$, where the $o_P(1)$ term is uniform for $1/(k+1)\leq
s\leq k/(k+1)$. The functions $a_0(m)$ and $A_0(m)$ are chosen as in
Lemma~\ref{2ndordlem} below.
%-\sqrt k A(m)\left( H_{\gamma,\rho} \left(\frac1{-\log s}\right)-
%for some functions $r_{k,m}$ and $R_{k,m}$ (the latter random), with $k
%as $n\to\infty$. Both functions $r_{k,m}$ and $R_{k,m}$ are integrable
%over [0,1].
\end{theo}

%th2.2 #&#
\begin{theo}\label{AsymptNMr}
Assume the conditions of Theorem~\ref{sup_lem} with $\gamma<1/2$. Then
\begin{eqnarray*}
&&\sqrt k \biggl(\frac{(r+1)\beta_r-b_m}{a_m}-D_r(\gamma) \biggr)
\\
&&\qquad\to^d (r+1)\int_0^1
s^{r-1}(-\log s)^{-1-\gamma}E(s) \,ds+\lambda I_r(\gamma,
\rho)=:Q_r,
\end{eqnarray*}
as $n\to\infty$, jointly for $r=0,1,2,3,\ldots,$ where $\to^d$ means
convergence in distribution, $E$ is Brownian bridge,
\[
D_r(\xi)=\frac{(r+1)^{\xi} \Gamma(1-\xi)-1}{\xi},\qquad \xi<1
\]
[$D_r(0)=\log(r+1)-\Gamma'(1)$ as defined by continuity], and
\[
I_r(\gamma,\rho)
=\cases{
\displaystyle\frac{1} \rho \bigl(D_r(
\gamma+\rho)-D_r(\gamma) \bigr), \vspace*{2pt}\cr
\qquad\rho\neq0,
\vspace*{2pt}\cr
\displaystyle D'_r(\gamma)
=\frac{(r+1)^\gamma}\gamma \bigl(-\Gamma'(1-\gamma) + \log(r+1)
\Gamma (1-\gamma)\vspace*{2pt}\cr
\hspace*{150pt}{}-(r+1)^{-\gamma}D_r(\gamma) \bigr),\vspace*{2pt}\cr
\qquad\gamma
\neq0,\rho= 0,
\vspace*{2pt}\cr
\displaystyle D'_r(0)=\frac{1} 2 \bigl(\log^2(r+1)+\Gamma''(1)-2
\log(r+1)\Gamma'(1) \bigr),\vspace*{2pt}\cr
\qquad \gamma= 0,\rho= 0.}
\]
\end{theo}

Note that $\Gamma'(1-\gamma)=\int_0^\infty u^{-\gamma}e^{-u}(\log u)
\,du$ and $\Gamma''(1)=1.97811$.\vadjust{\goodbreak}

%re2.1 #&#
\begin{rem} The condition $\sqrt k A(m)\to\lambda\in\R$ means that the
growth of $k_n$, the number of blocks, is restricted with respect to
the growth of $m_n$, the size of a block, as $n\to\infty$. In
particular this condition implies that $(\log k)/m\to0$, as $n\to\infty$.
\end{rem}

%The asymptotic normality of $\hat\gamma_{k,m}$, $\hat a_{k,m}$ and $

%th2.3 #&#
\begin{theo}\label{AsymptNgab}
Under the conditions of Theorem~\ref{AsymptNMr}, as $n\to\infty$,
\begin{eqnarray*}
\sqrt k (\hat\gamma_{k,m}-\gamma )&\to^d&\frac{1}{\Gamma
(1-\gamma)}
\biggl(\frac{\log3}{1-3^{-\gamma}}-\frac{\log2}{1-2^{-\gamma
}} \biggr)^{-1}
\\
&&{}\times \biggl\{\frac{\gamma}{3^{\gamma}-1} (Q_2-Q_0 )-
\frac{\gamma
}{2^{\gamma}-1} (Q_1-Q_0 ) \biggr\}\\
&=:&\Delta,
\\
%&&+\frac\lambda{\Gamma(1-\gamma)}\left(\frac{\log3}{1-3^{-\gamma}}-
%&&        \left\{\frac{\gamma}{3^{\gamma}-1}\left(I_2(\gamma,
\sqrt k \biggl(\frac{\hat a_{k,m}}{a_m}-1 \biggr)&\to^d&\frac\gamma {
\bigl(2^\gamma-1 \bigr)\Gamma(1-\gamma)} (Q_1-Q_0
)
\\
&&\hspace*{-4pt}{}+\Delta \biggl\{\frac{\log2}{\gamma} \biggl(\frac{-\gamma}{1-2^{-\gamma
}}+
\frac{1}{\log2} \biggr) +\frac{\Gamma'(1-\gamma)}{\Gamma(1-\gamma
)} \biggr\}\\
&=:&\Lambda,
\\
\sqrt k\frac{\hat b_{k,m}-b_m}{a_m}&\to^d& Q_0+
\frac{\gamma\Gamma'(1-\gamma)-1+\Gamma(1-\gamma)}{\gamma^2} \Delta+
\frac{1-\Gamma(1-\gamma)}{\gamma} \Lambda\\
&=:&
\Xi;
\end{eqnarray*}
where for $\gamma=0$ the formulas should read as (defined by continuity):
\begin{eqnarray*}
\sqrt k\hat\gamma_{k,m}&\to^d& \biggl(\frac{\log3}2-
\frac{\log2}2 \biggr)^{-1} \biggl(\frac{1}{\log3}
(Q_2-Q_0 )-\frac{1}{\log2} (Q_1-Q_0
) \biggr),
\\
\sqrt k \biggl(\frac{\hat a_{k,m}}{a_m}-1 \biggr) &\to^d&
\frac{1}{\log2} (Q_1-Q_0 )+\Delta \biggl(
\frac{\log
2}2+\Gamma'(1) \biggr),
\\
\sqrt k\frac{\hat b_{k,m}-b_m}{a_m} &\to^d& Q_0-
\Gamma''(1)\Delta+\Gamma'(1)\Lambda.
\end{eqnarray*}
\end{theo}

%re2.2 #&#
\begin{rem}
A slight modification of $\hat\gamma_{k,m}$ produces the explicit estimator
%
%e5 #&#
\begin{equation}
\hat\gamma_{k,m}^*=\frac{1}{\log2} \log \biggl(\frac{4\beta_3-\beta
_0}{2\beta_1-\beta_0}-1
\biggr),
\end{equation}
which is the solution of $
 (4^{\hat\gamma_{k,m}^*}-1 ) (2^{\hat\gamma
_{k,m}^*}-1 )^{-1}= (4\beta_3-\beta_0 ) (2\beta
_1-\beta_0 )^{-1}$.
The conditions of Theorem~\ref{AsymptNMr} imply
\begin{eqnarray*}
\sqrt k \bigl(\hat\gamma_{k,m}^*-\gamma \bigr)&\to^d&
\frac{1}{\Gamma(1-\gamma)} \biggl(\frac{\log4}{1-4^{-\gamma}}-\frac{\log
2}{1-2^{-\gamma}}
\biggr)^{-1}\\
&&{}\times \biggl\{\frac{\gamma}{4^{\gamma}-1} (Q_3-Q_0
)-\frac{\gamma
}{2^{\gamma}-1} (Q_1-Q_0 ) \biggr\}.
\end{eqnarray*}
\end{rem}

%s2.2 #&#
\subsection{High quantile estimation}
\label{Asymptnormq_sect}

Our estimator for $x_n=F^\leftarrow (1-p_n )= V (1/
(-\log (1-p_n ) ) )$, with $p_n$ small, is
\[
\hat x_{k,m}=\hat b_{k,m}+\hat a_{k,m}
\frac{ (mp_n )^{-\hat
\gamma_{k,m}}-1}{\hat\gamma_{k,m}}.
\]

%th2.4 #&#
\begin{theo}\label{AsymptNquantile}
Assume the conditions of Theorem~\ref{AsymptNMr} with $\rho$ negative,
or zero with $\gamma$ negative. Moreover, assume that the probabilities
$p_n$ satisfy
\[
\lim_{n\to\infty}mp_n=0 \quad\mbox{and}\quad \lim
_{n\to\infty}\frac
{\log(mp_n)}{\sqrt k}=0
\]
[in case $\rho<0$ the latter can be simplified to $\lim_{n\to\infty
}(\log p_n)/\sqrt k=0$]. Then
\[
\sqrt k\frac{ (\hat x_{k,m}-x_n )}{a_m q_\gamma
(1/(mp_n) )}\to^d \Delta+(\gamma_-)^2 \Xi-
\gamma_-\Lambda -\lambda\frac{\gamma_-}{\gamma_-+\rho}
\]
as $n\to\infty$, where $\gamma_-=\min(0,\gamma)$ and $q_\gamma(t)=\int_1^t s^{\gamma-1}\log s \,ds$.
\end{theo}

%distribution $F$. One may also want to estimate a high quantile of the
%distribution of the block maximum. In that case we need to estimate
%$x_n:=V\left(m/\left(-\log\left(1-p_n\right)\right)\right)$. The
%result is as above with
%$q_n^*:=-\log(1-p_n)/m$ replacing $q_n$ and $c_n^*:=1/mq_n^*$ replacing
%$c_n$ replaced by $mc_n$.

%s3 #&#
\section{Theoretical comparison between BM and POT methods}
\label{blockmaxPOT_sect}

In this section, we develop a theoretical comparison between the BM and
POT methods, by comparing the two PWM estimators for the two methods
[Hosking and Wallis (\citeyear{HoskingWallis87}) and
Hosking, Wallis and Wood (\citeyear{HoskingWallisWood85}), resp., for POT and BM].

First, we introduce the PWM-POT estimators for $\gamma$ and $a(n/k)$,
where $k$ is the number of selected order statistics, $\{\tilde
X_{n-i,n}\}_{i=0}^{k-1}$, from the original sample $\tilde X_1,\tilde
X_2,\ldots,\tilde X_n$. The statistics
\[
%begin{eqnarray*}
P_n= \frac{1}{k} \sum
_{i=0}^{k-1} \tilde X_{n-i,n} - \tilde
X_{n-k,n} \quad\mbox{and}\quad Q_n= \frac{1}{k} \sum
_{i=0}^{k-1} \frac{i}{k} (\tilde
X_{n-i,n} - \tilde X_{n-k,n} )
\]
%
%end{eqnarray*}
are estimators for $a(n/k)(1-\gamma)^{-1}$ and $a(n/k) (2(2-\gamma
)^{-1} )$, respectively. Consequently, the PWM estimators are
\[
%begin{eqnarray*}
\hat\gamma_{k,n}=1- \biggl(\frac{P_n}{2 Q_n}-1
\biggr)^{-1} \quad\mbox{and}\quad \hat a(n/k)=P_n \biggl(
\frac{P_n}{2 Q_n}-1 \biggr)^{-1}.
\]
%
%end{eqnarray*}
The quantile estimator is
\[
\hat x_{k,n}=\tilde X_{n-k,n}+\hat a(n/k)\frac{ ({k}/{(np_n)}
)^{\hat\gamma_{k,n}}-1}{\hat\gamma_{k,n}}.
\]
Asymptotic normality under conditions equivalent to the ones in
Theorems \ref{AsymptNgab} and~\ref{AsymptNquantile} holds [see, e.g.,
\citet{CaiHaanZhou}], if $\rho\in[-1,0]$ with a caveat for $\rho
=-1$ [for certain cases the functions $A$ in the corresponding
second-order conditions may not be the same asymptotically resulting in
different values of $\lambda$ in the limiting distributions;
cf. \citet{DreesHaanLi03}].

For BM, $k$ is defined as the number of blocks and, for POT, $k$ is
defined as the number of selected top order statistics. Hence, in both
cases $k$ means the number of selected observations. For the
theoretical comparison, we confine ourselves to the range $\rho\in
[-1,0]$ and $\gamma\in[-1,1/2)$, a usual range in many applications.

%In the POT approach, asymptotic normality of the PWM estimators are
%well known see e.g. \citet{CaiHaanZhou}. The conditions are
%basically the same in this case but with a second order condition
%formulated in terms of the function $U=(1/(1-F))^\leftarrow$, instead
%of the function $V=\left(-1/\log F\right)^{\leftarrow}$. For a
%comparison between the two second order conditions see Drees, de Haan
%and Li (2003). In particular for $\rho\in[-1,0]$ (with a caveat for $
%second order parameter $\rho$ and auxiliary function. Hence, the
%conclusions below on comparing the asymptotic properties of the
%estimators are valid within this range of $\rho$. We also confine to
%the range $\gamma\in[-1,1/2)$ which seems the most usual range in
%applications. Recall that asymptotic normality is only valid for $

\subsection*{Extreme value index estimators}

%of $k$ which is different in the two cases. The `optimal choice' is
%the value that makes the limiting mean square error of $\hat\gamma-
%and can choose the $k$ at will. We also compare the two optimal
%$k_0$'s.

%In Theorem~\ref{AsymptNgab} the asymptotic normality of the PWM
%estimator under the BM approach was established.

%As explained, under basically the same conditions of Theorem~\ref{AsymptNgab} we have for the PWM estimators under POT,
%as $n\to\infty$, with $N$ normal random variable and $\lambda$ as in
%Theorem~\ref{AsymptNgab}.

%
 \begin{itemize}
\item First, we compare asymptotic variance and bias for a common
value of $k$:

The asymptotic variances of the two $\gamma$ estimators are shown in
Figure~\ref{VarBMPOT.fig}: the curve from BM is always below the other
one, meaning lower values for the asymptotic variance for all values of
$\gamma$. The asymptotic biases are compared in Figure~\ref{biasRATIO.fig}, through the ratio
``bias BM/bias POT''. Recall that the
bias depends on both first- and second-order parameters $\gamma$ and
$\rho$. Contrary to what is observed for the variance, the bias of BM
is always larger but for $\rho=0$ they are the same regardless the
value of $\gamma$, equal to 1 [or $\lambda$ if one takes into account
the asymptotic contribution of $\sqrt k A(n/k)$ to the biases]. %In
%Figures \ref{biasBM.fig} and \ref{biasPOT.fig} are also represented
%the asymptotic bias for each case separately.

%f1 #&#
\begin{figure}

\includegraphics{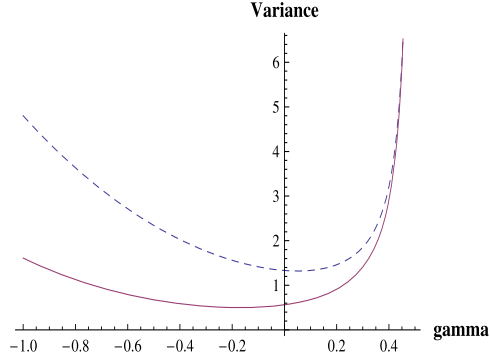}

\caption{Asymptotic variances of $\gamma$ PWM estimators with dashed
line for POT.} \label{VarBMPOT.fig}
\end{figure}

%f2 #&#
\begin{figure}[b]

\includegraphics{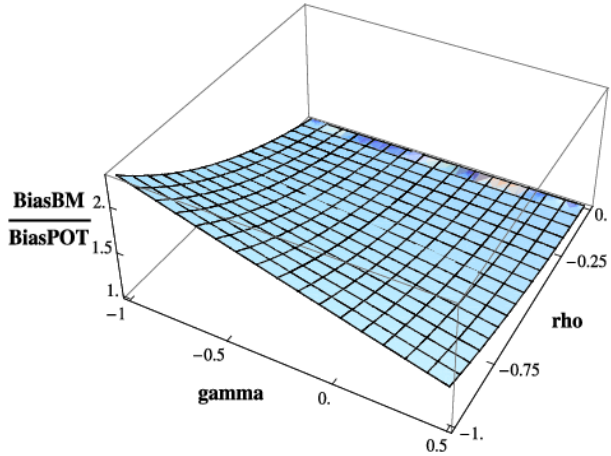}

\caption{Ratio of asymptotic bias of $\gamma$ PWM estimators.} \label
{biasRATIO.fig}
\end{figure}

%f3 #&#
\begin{figure}

\includegraphics{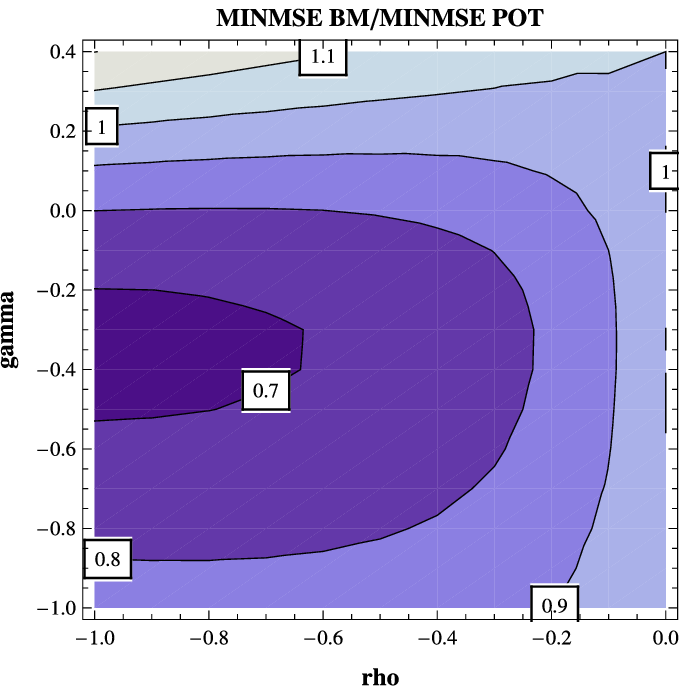}

\caption{Contour plot for the ratio of asymptotic minimal mean square
error of $\gamma$ PWM estimators.} \label{ratiocontour.fig}
\end{figure}

\item Next, we compare asymptotic mean square errors for
the ``optimal choice'' of $k$ (i.e., that value that makes the limiting
mean square error of $\hat\gamma-\gamma$ minimal), which is different
in the two cases:

An asymptotic expression of the ``asymptotic minimal mean square error''
(MINMSE in the sequel) is obtained in the following way.
Suppose $\rho<0$. First we find for each estimator the optimal $k$ in
the sense of minimizing the approximate asymptotic mean square error.
Denote by $\sigma^2_i=\sigma^2_i(\gamma)$ and $B^2_i=B^2_i(\gamma,\rho
)$ ($i=1,2$; ``1'' refers to PWM-BM and ``2'' refers to PWM-POT) the
asymptotic variance and squared bias of the estimators. %and recall the
%$\rho$ regularly varying function $A(\cdot)$ from the second order
Under Condition \ref{2ndordcond}, we can write $A^2(t)=\int_t^\infty
s(u) \,du$ with $s(\cdot)$ decreasing and $2\rho-1$ regularly varying.
The limiting mean square error is, approximately,
%
%e6 #&#
\begin{equation}
\label{minmse} \inf_k \biggl(\frac{\sigma_i^2}k+A^2(n/k)B_i^2
\biggr)
\end{equation}
or, writing $r$ for $n/k$, $ \inf_r ( (r/n)\sigma_i^2+B_i^2\int_r^\infty s(u) \,du )$.
Setting the derivative equal to zero and using properties of regularly
varying functions one finds for the optimal choice of $r$,
$ r_0^{(i)}\sim (1/ s )^\leftarrow(n)   (B_i^2/\sigma
_i^2 )^{1/(1-2\rho)} $
and, in terms of $k$,
\[
k_0^{(i)}\sim\frac{n}{ ({1} / s )^{\leftarrow}(n)} \biggl(\frac{\sigma_i^2}{B_i^2}
\biggr)^{1/(1-2\rho)}.
\]
Note that the optimal $k_0^{(i)}$ is different but of the same order
for both methods. Next, inserting $k_0^{(i)}$ in \eqref{minmse}, after
some manipulation we get the following asymptotic expression for MINMSE,
\[
\frac{1-2\rho}{-2\rho}\frac{ ({1}/ s )^{\leftarrow}(n)}n \bigl(B_i^2
\bigr)^{1/(1-2\rho)} \bigl(\sigma_i^2
\bigr)^{-2\rho/(1-2\rho)}.
\]
It follows that MINMSE(BM)/MINMSE(POT) is, approximately,
\[
 \biggl(\frac{B_1^2(\gamma,\rho)}{B_2^2(\gamma,\rho)} \biggr)^{1/(1-2\rho
)} \biggl(
\frac{\sigma_1^2(\gamma)}{\sigma_2^2(\gamma)} \biggr)^{-2\rho
/(1-2\rho)},
\]
which does not depend on n, just on $\gamma$ and $\rho$.

%The ratio `MINMSE(BM)/MINMSE(POT)' is represented in Figure~\ref{ratioMINMSE.fig} (??Do we keep it??) and the corresponding
The contour plot of ``MINMSE(BM)/MINMSE(POT)'' is represented in Figure~\ref{ratiocontour.fig}. It can be seen that the BM has lower MINMSE for
a large range of $(\gamma,\rho)$ combinations. Note that this range
includes $\gamma$ negative and $\gamma$ positive close to zero which
seem to be common values in many practical situations, for example, in
hydrologic and climatologic extremes. Only for $\gamma>0.2$
approximately, MINMSE for POT can be lower depending on $\rho$.

%estimators.} \label{ratioMINMSE.fig}

%f4 #&#
\begin{figure}

\includegraphics{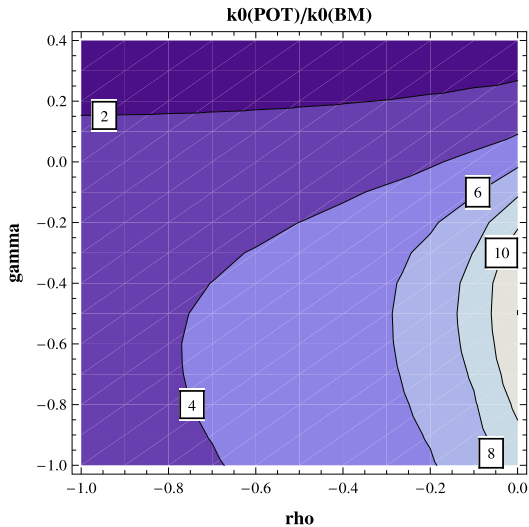}

\caption{Contour plot for the ratio of the optimal values of $k$.}
\label{contourk0RATIO.fig}
\end{figure}

%The last picture, Figure~\ref{contourk0RATIO.fig} corresponds to the
%ratio of the optimal values of $k$.
Finally, comparing the optimal sample sizes (cf. Figure~\ref{contourk0RATIO.fig} with contour plot of the ratio of the optimal
values of $k$), one sees that POT requires systematically larger
optimal sample size even when the approximate MINMSE is smaller for POT
than BM.
\end{itemize}

\subsection*{Quantile estimators}

We repeat the previous analysis for the quantile estimators:
%In Theorem~\ref{AsymptNquantile} the asymptotic normality of the PWM
%quantile estimator under the BM approach was established.

%In the POT approach, asymptotic normality of the PWM quantile
%estimators are well known see e.g. \citet{CaiHaanZhou}.
%Similarly as before, under basically the same conditions of Theorem~\ref{AsymptNquantile} we have for the PWM quantile estimators under
%POT,
%+\frac{(2-\gamma)^2(1-\gamma+\gamma_-)^2(4+3\gamma^2-2\gamma^3)}{(3-2
%as $n\to\infty$.

\begin{figure}[b]

\includegraphics{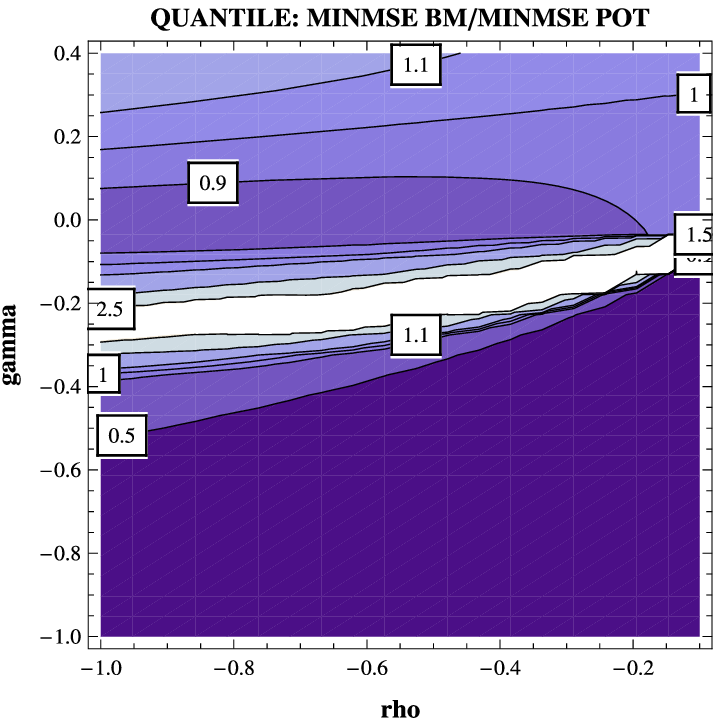}

\caption{Contour plot for the ratio of asymptotic minimal mean square
error of quantile PWM estimators.} \label{contourQRATIO.fig}
\end{figure}

\begin{figure}

\includegraphics{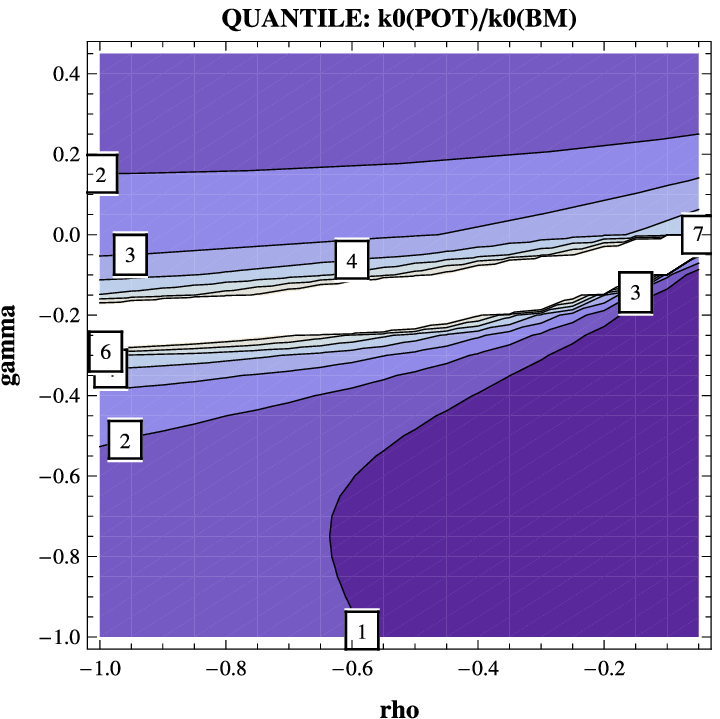}

\caption{Contour plot for the ratio of the optimal values of $k$ in
quantile case.} \label{contourk0QRATIO.fig}
\end{figure}

%Note that the asymptotic distributions are the same as in $\gamma$
%estimation when $\gamma\geq0$.
\begin{itemize}
\item The asymptotic variances of the two estimators are
compared in Figure~\ref{VarQBMPOT.fig}: again the curve from BM is
always below the other one meaning lower values for the asymptotic
variance for all values of $\gamma$. In Figures~\ref{biasQBM.fig} and
\ref{biasQPOT.fig}, the asymptotic bias is represented for each case
separately. Note that for $\gamma$ negative, the bias for BM approaches
zero when $\rho\uparrow0$ whereas in the POT case it escapes to
$-\infty$.
%f5 #&#
\begin{figure}

\includegraphics{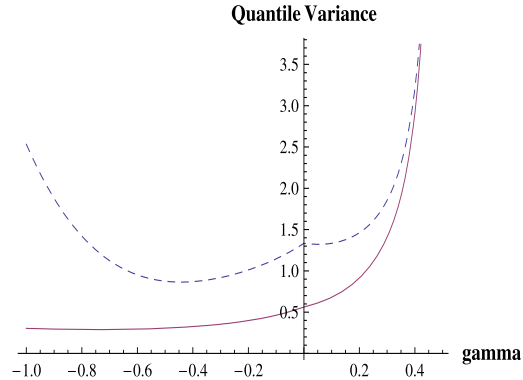}

\caption{Asymptotic variances of quantile PWM estimators with dashed
line for POT.} \label{VarQBMPOT.fig}
\end{figure}

%f6 #&#
\begin{figure}[b]

\includegraphics{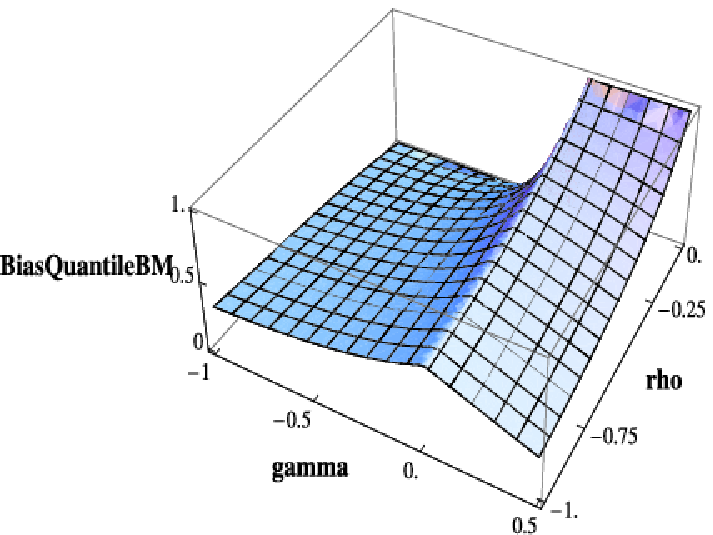}

\caption{Asymptotic bias of quantile PWM estimator under BM method.}
\label{biasQBM.fig}
\end{figure}

%f7 #&#
\begin{figure}

\includegraphics{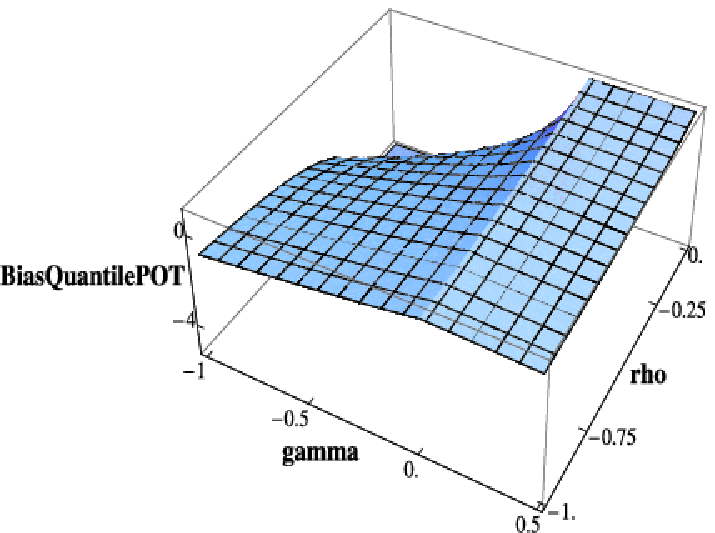}

\caption{Asymptotic bias of quantile PWM estimator under POT method.}
\label{biasQPOT.fig}
\end{figure}

\item The contour plot for the ratio
``MINMSE(BM)/MINMSE(POT)'' is represented in Figure~\ref{contourQRATIO.fig}. Again the BM method has lower MINMSE for a large
range of $(\gamma,\rho)$ combinations. The ``irregularity'' around $\gamma
\approx-0.2$ is due to a change of sign in the bias in the POT case.
Finally, Figure~\ref{contourk0QRATIO.fig} gives the contour plot for
the ratio of the optimal values of $k$, which is smaller than one when
$\gamma$ is small and $\rho$ is closer to zero.
\end{itemize}

%f8 #&#

%f9 #&#

\textit{In conclusion}, for both the extreme value index and quantile PWM
estimators, the ones from the BM method have always lower asymptotic
variances. Moreover, at an optimal level the BM gives lower MINMSE,
thus being more efficient, under many practical situations. This is in
agreement with some of Sofia Caires' (\citeyear{Caires09}) conclusions, for example,
that for equal sample sizes or with more than 200 years of data the
uncertainty or the error of the estimates are lower for BM than for POT.

%s4 #&#
\section{Proofs}
\label{proofs_sect}
Throughout this section, $Z$ represents a unit Fr\'echet random
variable, that is, one with distribution function $F(x)=e^{-1/x}$,
$x>0$, and $ \{ Z_{i,k} \}_{i=1}^k$ are the order statistics
% (define $Z_{0,k}:=Z_{1,k}$)
from the associated i.i.d. sample of size $k$, $Z_1,\ldots,Z_k$.
Similarly, $ \{ X_{i,k} \}_{i=1}^k$
%(again define $X_{0,k}:=X_{1,k}$)
represents the order statistics of the block maxima $X_1,\ldots,X_k$
from \eqref{blockmaxima_def} and, $X_{\lceil u\rceil,k}=X_{r,k}$ for
$r-1<u\leq r$, $r=1,\ldots,k$.
%represents the smallest integer larger or equal to $x\in\R$)
Recall the function $V$ from Section~\ref{AsympN_sect}. The following
representation will be useful:
%
%e7 #&#
\begin{equation}
\label{XVZrepresentation} X=^d V(mZ).
\end{equation}

We start by formulating a number of auxiliary results.

%le4.1 #&#
\begin{lem}\label{Zaux_lem}
1. As $k\to\infty$,
\[
(\log k)Z_{1,k}\to^P 1.\vspace*{-9pt}% \mbox{ and }
\]
\begin{enumerate}[2.]
\item[2.][Cs\"org\H{o} and Horv\'ath (\citeyear{Csorgo}), page~381]
Let $0<\nu<1/2$. With
$\{E_k\}_{k\geq1}$, an appropriate sequence of Brownian bridges,
\[
\sup_{1/(k+1)\leq s\leq k/(k+1)}\frac{s(-\log s)}{ (s(1-s)
)^\nu} \biggl\llvert \sqrt{k} \bigl( (-
\log s) Z_{\lceil ks\rceil,k}-1 \bigr)-\frac
{E_k(s)}{s(-\log s)}\biggr\rrvert =
o_P(1),
\]
as $k\to\infty$ ($\lceil u\rceil$ represents the smallest integer
larger or equal to $u$).
\item[3.] Similarly, with $0<\nu<1/2$ for an appropriate
sequence $\{E_k\}_{k\geq1}$ of Brownian bridges and $\xi\in\R$,
\begin{eqnarray*}
&&\sup_{1/(k+1)\leq s\leq k/(k+1)} \bigl(s(1-s) \bigr)^{-\nu}
\\
&&\hspace*{-4pt}\qquad{}\times\biggl| \sqrt{k} s(-\log s)^{1+\xi} \biggl( \frac{Z_{\lceil ks\rceil,k}^\xi
-1}{\xi}-
\frac{(-\log s)^{-\xi}-1}{\xi} \biggr) -E_k(s) \biggr| = o_P(1),% \mbox{ as } k\to\infty.
\end{eqnarray*}
as $k\to\infty$.
\end{enumerate}
\end{lem}

The following is an easily obtained variant of Theorem B.3.10 of \citet{HaanFerreira}.

%le4.2 #&#
\begin{lem}\label{2ndordlem}
Under Condition \ref{2ndordcond}, there are functions $A_0(t)\sim A(t)$
and $a_0(t)=a(t) (1+o (A(t) ) )$, as $t\to\infty$,
such that for all $\varepsilon,\delta>0$ there exists
$t_0=t_0(\varepsilon,\delta)$ such that for $t,tx> t_0$,
%
%e8 #&#
\begin{eqnarray}
\label{2ndordcondunif} &&\biggl\llvert \frac{ {(V(tx)-V(t))}/{a_0(t)}-
{(x^\gamma-1)}/{\gamma}
}{A_0(t)}-H_{\gamma,\rho}(x)\biggr\rrvert
\nonumber
\\[-8pt]
\\[-8pt]
\nonumber
&&\qquad
\leq\varepsilon\max \bigl( x^{\gamma+\rho+\delta},x^{\gamma+\rho-\delta} \bigr).
\end{eqnarray}
Moreover,
%
%e9 #&#
\begin{eqnarray}
\label{2ndordcondunifa} &&\biggl\llvert \frac{ {a_0(tx)}/{a_0(t)} -
x^\gamma}{A_0(t)}- x^\gamma
{(x^\rho-1)}/{\rho} \biggr\rrvert
\nonumber
\\[-8pt]
\\[-8pt]
\nonumber
&&\qquad\leq\varepsilon\max \bigl(x^{\gamma+\rho
+\delta},x^{\gamma+\rho-\delta}
\bigr)
\end{eqnarray}
and
\[
\biggl\llvert \frac{A_0(tx)}{A_0(t)} - x^\rho\biggr\rrvert \leq\varepsilon
\max \bigl(x^{\rho+\delta},x^{\rho-\delta}\bigr).
\]
\end{lem}

Note that
\[
H_{\gamma,\rho}(x)=\cases{ %
\displaystyle\frac{1}{\rho} \biggl(
\frac{x^{\gamma+\rho}-1}{\gamma+\rho} - \frac
{x^{\gamma}-1}{\gamma} \biggr) , & \quad $\rho\neq0\neq\gamma,$
\vspace*{2pt}\cr
\displaystyle\frac{1}{\gamma} \biggl(x^\gamma\log x - \frac{x^{\gamma}-1}{\gamma} \biggr) ,
&\quad  $\rho=0\neq\gamma,$
\vspace*{2pt}\cr
\displaystyle\frac{1}{\rho} \biggl(\frac{x^{\rho}-1}{\rho}-\log x \biggr) , &\quad
$\rho\neq0=\gamma,$
\vspace*{2pt}\cr
\displaystyle\frac{1}{2} (\log x)^2 , &\quad $\rho= 0=\gamma.$}
\]

\begin{pf*}{Proof of Theorem~\ref{sup_lem}}
By representation \eqref{XVZrepresentation},
\begin{eqnarray*}
&&\frac{X_{\lceil ks\rceil,k}-b_m}{a_0(m)}- \frac{(-\log
s)^{-\gamma}-1}{\gamma}
\\
&&\qquad =^d \biggl(\frac{V (mZ_{\lceil ks\rceil,k} )-b_m}{a_0(m)}-\frac
{V ({m}/{-\log s} )-b_m}{a_0(m)} \biggr) \\
&&\qquad\quad{}+ \biggl(
\frac{V ({m}/{-\log s} )-b_m}{a_0(m)}-\frac{(-\log
s)^{-\gamma}-1}{\gamma} \biggr)
\\
&&\qquad= \mbox{I (random part)} + \mbox{II (bias part)}.
\end{eqnarray*}

We start with part I,
\begin{eqnarray*}
\mbox{I}&=& \biggl\{ (-\log s)^{-\gamma}
\frac{V ((-\log s)Z_{\lceil ks\rceil,k}{m}/{-\log s} )-V ({m}/{-\log s} )}
{a_0 ({m}/{-\log s} )} \biggr\} \\
&&{}\times\biggl
\{\frac{a_0 ({m}/{-\log s} )}{a_0(m)}(-\log s)^{\gamma} \biggr\}
\\
&=&\mbox{I.1} \times\mbox{I.2}.
\end{eqnarray*}
According to \eqref{2ndordcondunifa} of Lemma~\ref{2ndordlem}, for each
$\varepsilon,\delta>0$ there exists $t_0$ such that the factor~I.2 is
bounded (above and below) by
\[
1+A_0(m) \biggl\{\frac{(-\log s)^{-\rho}-1}{\rho}\pm\varepsilon\max \bigl((-\log
s)^{-\rho+\delta},(-\log s)^{-\rho-\delta} \bigr) \biggr\}
\]
provided $m\geq t_0$ and $s\geq e^{-m/t_0}$. According to \eqref
{2ndordcondunif} of Lemma~\ref{2ndordlem}, for factor I.1 we have the bounds
\begin{eqnarray*}
&&(-\log s)^{-\gamma}\frac{ ((-\log s)Z_{\lceil ks\rceil,k}
)^{\gamma}-1}{\gamma}+A_0 \biggl(
\frac{m}{-\log s} \biggr) (-\log s)^{-\gamma
}
\\
&&\quad{}\times \bigl\{ H_{\gamma,\rho} \bigl((-\log s)Z_{\lceil ks\rceil,k} \bigr) \pm
\varepsilon\max \bigl( \bigl((-\log s)Z_{\lceil ks\rceil,k} \bigr)^{\gamma+\rho+\delta},\\
&&\hspace*{174pt}\bigl((-\log s)Z_{\lceil ks\rceil,k} \bigr)^{\gamma+\rho-\delta} \bigr) \bigr\}
\\
&&\qquad= \mbox{I.1a} + \mbox{I.1b}
\end{eqnarray*}
provided $s\geq e^{-m/t_0}$ and $m/\log k\geq t_0$ [the latter
inequality eventually holds true since $\sqrt k A_0(m)$ is bounded].
Note that $m/\log k\geq t_0$ implies $mZ_{1,k}\geq2t_0$ which implies
(Lemma~\ref{Zaux_lem}) $mZ_{\lceil ks\rceil,k}\geq2t_0$ for all $s$.

For term I.1a, we use Lemma~\ref{Zaux_lem}.3:
\[
\frac{ (Z_{\lceil ks\rceil,k} )^{\gamma}-1}{\gamma}-\frac
{(-\log s)^{-\gamma}-1}{\gamma}
\]
is bounded (above and below) by
\[
\frac{1}{\sqrt k} \frac{E_k(s)}{s(-\log s)^{1+\gamma}}\pm\frac
{\varepsilon}{\sqrt k}\frac{ (s(1-s) )^\nu}{s(-\log
s)^{1+\gamma}},
\]
for some $\varepsilon>0$, $0<\nu<1/2$ and all $s\in[1/(k+1),k/(k+1)]$.

Next, we turn to term I.1b. By Lemma~\ref{2ndordlem}, $(-\log
s)^{-\gamma}A_0 (\frac{m}{-\log s} )$ is bounded (above and
below) by
\[
A_0(m) \bigl\{(-\log s)^{-\gamma-\rho} \pm\varepsilon\max \bigl((-
\log s)^{-\gamma-\rho+\delta},(-\log s)^{-\gamma-\rho-\delta} \bigr) \bigr\}\vadjust{\goodbreak}
\]
provided $s>e^{-m/t_0}$ and $m/\log k>t_0$. Furthermore for $\rho\neq0\neq\gamma$ and $s\in[1/(k+1),k/(k+1)]$, by Lemma~\ref
{Zaux_lem}.3,
\begin{eqnarray*}
&&H_{\gamma,\rho} \bigl((-\log s)Z_{\lceil ks\rceil,k}\bigr)
\\
&&\qquad=\frac{1} \rho\biggl\{ \frac{((-\log s)Z_{\lceil ks\rceil,k})^{\gamma+\rho
}-1}{\gamma+\rho} -\frac{((-\log s)Z_{\lceil ks\rceil,k}
)^{\gamma}-1}{\gamma} \biggr\}
\\
&&\qquad=\frac{1} \rho\biggl\{(-\log s)^{\gamma+\rho}\biggl[
\frac{Z_{\lceil ks\rceil,k}^{\gamma+\rho}-1}{\gamma+\rho}-\frac{(-\log
s)^{-\gamma-\rho}-1}{\gamma+\rho}\biggr]
\\
&&\hspace*{24pt}\qquad\quad{} -(-\log s)^{\gamma}\biggl[ \frac{Z_{\lceil ks\rceil,k}^{\gamma}-1}{\gamma}-\frac{(-\log s)^{-\gamma
}-1}{\gamma}\biggr]
\biggr\}
\end{eqnarray*}
is bounded by
\begin{eqnarray*}
&&\frac{1} \rho\biggl\{(-\log s)^{\gamma+\rho} \biggl[ \frac{1}{\sqrt k}
\frac{E_k(s)}{s(-\log s)^{1+\gamma+\rho}}\pm\frac
{\varepsilon}{\sqrt k}\frac{ (s(1-s) )^\nu}{s(-\log
s)^{1+\gamma+\rho}} \biggr]
\\
&&\hspace*{24pt}\quad{}- (-\log s)^{\gamma} \biggl[ \frac{1}{\sqrt k} \frac{E_k(s)}{s(-\log s)^{1+\gamma}}\mp
\frac
{\varepsilon}{\sqrt k}\frac{ (s(1-s) )^\nu}{s(-\log
s)^{1+\gamma}} \biggr]\biggr\}
\\
&&\qquad=\pm\frac{2\varepsilon}{\rho\sqrt k}\frac{ (s(1-s) )^\nu
}{s(-\log s)},
\end{eqnarray*}
%
%Note that this expression is $O_P(1/\sqrt k)$ for fixed $s$.
and similarly for cases other than $\rho\neq0\neq\gamma$. The
remaining part of I.1b, namely $\pm\varepsilon\max ( ((-\log
s)Z_{\lceil ks\rceil,k} )^{\gamma+\rho+\delta}, ((-\log
s)Z_{\lceil ks\rceil,k} )^{\gamma+\rho-\delta} )$, is similar.
%and is at most
%0_P(1) \frac{B_k(s)}{\sqrt k s(-\log s)}+\frac{o_P(1)}{\sqrt k}\frac{

Part II, by the inequalities of Lemma~\ref{2ndordlem}, is bounded by
\[
A_0(m) \biggl\{ H_{\gamma,\rho} \biggl(\frac{1}{-\log s} \biggr)
\pm\varepsilon\max \bigl((-\log s)^{-\gamma-\rho+\delta},
(-\log s)^{-\gamma-\rho-\delta} \bigr)
\biggr\}
\]
hence it contributes $\sqrt k A_0(m) H_{\gamma,\rho}  (\frac{1}{-\log s} )$ to the result.

Collecting all the terms, one finds the result.% where the functions
%$r_{k,m}$ and $R_{k,m}$ can be taken as
%r_{k,m}(s)=\frac{\left(s(1-s)\right)^\nu}{s(-\log s)^{\xi_1}e^{-
%s)|}
%and
%R_{k,m}(s)=\frac{c_2 B_k(s)}{s(-\log s)^{\xi_3}e^{-\delta_3|\log(-\log
%s)|}},
%for some reals $c_1,c_2,\xi_1,\xi_2,\xi_3,\delta_1,\delta_2,\delta_3$.
%In particular, these constants are such that both $r_{k,m}$ and
%$R_{k,m}$ are integrable over [0,1].
\end{pf*}

\begin{pf*}{Proof of Theorem~\ref{AsymptNMr}}
Let, for $r=0,1,2,3,\ldots,$
\[
J_k^{(r)}(s)=\frac{(\lceil ks\rceil-1)\cdots(\lceil ks\rceil
-r)}{(k-1)\cdots(k-r)},\qquad  s\in[0,1].
\]
Note that $J_k^{(r)}(s)\to s^r$, as $k\to\infty$, uniformly in $s\in
[0,1]$, and
\begin{eqnarray*}
\frac{1} k\sum_{i=1}^k
\frac{(i-1)\cdots(i-r)}{(k-1)\cdots(k-r)}&=&\int_0^1
J_k^{(r)}(s)\,ds=\frac{1}{r+1}\\
&=&\int_0^1
s^r \,ds.
\end{eqnarray*}
Then
\begin{eqnarray*}
&&\sqrt k \biggl(\frac{(r+1)\beta_r-b_m}{a_m}-\frac{(r+1)^\gamma
\Gamma(1-\gamma)-1}{\gamma} \biggr)
\\
&&\qquad=\sqrt k \biggl(\frac{(r+1)\int_0^1 X_{\lceil ks\rceil,k}J_k^{(r)}(s)
\,ds-b_m}{a_m}\\
&&\hspace*{55pt}{}-(r+1)\int_0^1
\frac{(-\log s)^{-\gamma}-1}\gamma s^r \,ds \biggr)
\\
&&\qquad=\sqrt k (r+1) \int_0^1 \biggl(
\frac{X_{\lceil ks\rceil
,k}-b_m}{a_m}-\frac{(-\log s)^{-\gamma}-1}\gamma \biggr) J_k^{(r)}(s)
\,ds
\\
&&\qquad\quad{}-\sqrt k (r+1) \int_0^1
\frac{(-\log s)^{-\gamma}-1}\gamma \bigl(s^r-J_k^{(r)}(s)
\bigr) \,ds
\\
&&\qquad=\sqrt k (r+1) \int_0^{1/(k+1)} \biggl(
\frac{X_{\lceil ks\rceil
,k}-b_m}{a_m}-\frac{(-\log s)^{-\gamma}-1}\gamma \biggr)J_k^{(r)}(s)
\,ds
\\
&&\qquad\quad{}+\sqrt k (r+1) \int_{1/(k+1)}^{k/(k+1)} \biggl(
\frac{X_{\lceil ks\rceil
,k}-b_m}{a_m}-\frac{(-\log s)^{-\gamma}-1}\gamma \biggr)J_k^{(r)}(s)
\,ds
\\
&&\qquad\quad{}+\sqrt k (r+1) \int_{k/(k+1)}^1 \biggl(
\frac{X_{\lceil ks\rceil
,k}-b_m}{a_m}-\frac{(-\log s)^{-\gamma}-1}\gamma \biggr)J_k^{(r)}(s)
\,ds
\\
&&\qquad\quad{}-\sqrt k (r+1) \int_0^1
\frac{(-\log s)^{-\gamma}-1}\gamma \bigl(s^r-J_k^{(r)}(s)
\bigr)\,ds
\\
&&\qquad=\mbox{I.1}+\mbox{I.2}+\mbox{I.3}+\mbox{I.4}.
\end{eqnarray*}

For I.4: since $ (s^r-J_k^{(r)}(s) )=O(1/k)$ uniformly in $s$,
$\mathrm{I.4}=O(1/\sqrt k)$.

For I.1, note that
%
%e10 #&#
\begin{equation}
\label{leftover} \int_0^{1/(k+1)} \sqrt k
\frac{X_{\lceil ks\rceil,k}-b_m}{a_m} s^r \,ds=o_P(1).
\end{equation}
This follows since, the left-hand side of \eqref{leftover} equals, in
distribution,
\[
\frac{\sqrt k}{(k+1)^{r+1}} \frac{V ( m Z_{1,k} )-V(m)}{a_m}
\]
which, by Lemmas \ref{Zaux_lem}.1, \ref
{2ndordlem} and the fact that $m/\log k\to\infty$, is bounded (below
and above) by
\[
\frac{\sqrt k}{(k+1)^{r+1}} \biggl\{\frac{Z_{1,k}^\gamma-1}{\gamma
}+A(m)H_{\gamma,\rho} (
Z_{1,k} )\pm\varepsilon A(m) \max \bigl(Z_{1,k}^{\gamma+\rho+\delta},Z_{1,k}^{\gamma+\rho-\delta}
\bigr) \biggr\}.
\]
This is easily seen to converge to zero in probability, since
$Z_{1,k}^{\xi}/\sqrt k = \break \{(\log k)  Z_{1,k}\}^{\xi}\log^{-\xi} k
/\sqrt k\to^P 0$ for all real $\xi$ and $\sqrt k A(m)\to\lambda$.
Hence, $\mathrm{I.1}=o_P(1)$.

Next, we show that
%
%e11 #&#
\begin{equation}
\label{leftover2} \int_{k/(k+1)}^1 \sqrt k
\frac{X_{\lceil ks\rceil,k}-b_m}{a_m} J_k^{(r)}(s) \,ds=o_P(1).
\end{equation}
The left-hand side equals, in distribution, since $J_k^{(r)}(s)\equiv
1$ for $s\in(k(k+1)^{-1},1)$,
\[
\biggl(1-\frac{k}{k+1} \biggr)\sqrt k \frac{V ( m Z_{k,k} )-V(m)}{a_m}.
\]
Lemma~\ref{Zaux_lem} yields
\[
\frac{V ( m Z_{k,k} )-V(m)}{a_m}=\frac{Z_{k,k}^\gamma-1}{\gamma
}+A(m) \bigl\{H_{\gamma,\rho} (
Z_{k,k} )\pm\varepsilon Z_{k,k}^{\gamma+\rho+\delta} \bigr\},
\]
which is (since $Z_{k,k}^\gamma/k^\gamma$ converges to a positive
random variable) of the order $O_P (k^\gamma )$. Hence, the
integral is of order $(k+1)^{-1}\sqrt k k^\gamma$ which tends to zero
since $\gamma<1/2$.

Finally, I.2 has the same asymptotic behaviour as
\[
(r+1)\int_{1/(k+1)}^{k/(k+1)} \sqrt k \biggl(
\frac{X_{\lceil ks\rceil
,k}-b_m}{a_m}-\frac{(-\log s)^{-\gamma}-1}{\gamma} \biggr) s^r \,ds,
\]
which, by Theorem~\ref{sup_lem} tends to
\[
(r+1)\int_0^1 s^{r-1}(-\log
s)^{-1-\gamma}E(s) \,ds+\lambda(r+1)\int_0^1
H_{\gamma,\rho} \biggl(\frac{1}{-\log s} \biggr) s^r \,ds.
\]
For the evaluation of the latter integral note that for $\xi<1$,
\[
(r+1)\int_0^1 s^r(-\log
s)^{-\xi}\,ds=(r+1)^{\xi-1}\int_0^\infty
v^{-\xi
}e^{-v}\,dv=(r+1)^{\xi-1}\Gamma(1-\xi).
\]
Moreover, note that
\[
(r+1)\int_0^1 s^r
\frac{(-\log s)^{-\xi}-1}\xi \,ds=\frac{(r+1)^{\xi}
\Gamma(1-\xi)-1}{\xi},\qquad \xi<1
\]
[$D_r(0)=\log(r+1)-\Gamma'(1)$ as defined by continuity], and
$(r+1)\times\break \int_0^1 H_{\gamma,\rho} (\frac{1}{-\log s}
)s^r \,ds=I_r(\gamma,\rho)$.
% \,ds=}\\
%&=&\left\{\begin{array}{ll}
% D'_r(\gamma)=& \\
% =\frac{(r+1)^\gamma}\gamma\left(-\Gamma'(1-\gamma) + \log(r+1)
%D'_r(0)=& \\
%=\frac1 2\left(\log^3(r+1)+\Gamma''(1)-3\log(r+1)\Gamma'(1)\right),&
\end{pf*}

\begin{pf*}{Proof of Theorem~\ref{AsymptNgab}}
From Theorem~\ref{AsymptNMr}, we obtain
\begin{eqnarray*}
\sqrt{k} \biggl(\frac{2\beta_1-\beta_0}{a_m}-\frac{2^\gamma-1}{\gamma
}\Gamma(1-\gamma) \biggr)
&\to^d &Q_1-Q_0,\\
\sqrt{k} \biggl(\frac{3\beta_2-\beta_0}{a_m}-\frac{3^\gamma-1}{\gamma
}\Gamma(1-\gamma) \biggr)
&\to^d& Q_2-Q_0
\end{eqnarray*}
hence, by Cram\'er's delta method,
\begin{eqnarray*}
&&\sqrt{k} \biggl(\frac{3^{\hat\gamma_{k,m}}-1}{2^{\hat\gamma
_{k,m}}-1}-\frac{3^{\gamma}-1}{2^{\gamma}-1} \biggr)\\
&&\qquad =\sqrt{k}
\biggl(\frac{3\beta_2-\beta_0}{2\beta_1-\beta_0}-\frac{3^\gamma
-1}{2^{\gamma}-1} \biggr)
\\
&&\qquad\to^d\frac{1}{\Gamma(1-\gamma)} \frac{3^\gamma-1}{2^{\gamma}-1} \biggl(\frac
\gamma{3^\gamma-1} (Q_2-Q_0 ) -\frac
\gamma{2^\gamma-1} (Q_1-Q_0 ) \biggr).
%&&+\frac\lambda{\Gamma(1-\gamma)} \frac{3^\gamma-1}{2^{\gamma}-1}
%-\frac\gamma{2^\gamma-1}\left(I_1(\gamma,\rho)-I_0(\gamma,\rho)\right)
\end{eqnarray*}
It follows that $\hat\gamma_{k,m}\to^P\gamma$, and hence
\[
\sqrt{k} \biggl(\frac{r^{\hat\gamma_{k,m}}-1}{r^{\gamma}-1}-1 \biggr)=\sqrt k \frac{r^{\hat\gamma_{k,m}-\gamma}-1}{1-r^{-\gamma}}
\]
has the same limit distribution as
\[
\sqrt{k} (\hat\gamma_{k,m}-\gamma )\frac{\log r}{1-r^{-\gamma
}},\qquad r=2,3.
\]
It follows that
\begin{eqnarray*}
&&\sqrt{k} \biggl(\frac{3^{\hat\gamma_{k,m}}-1}{2^{\hat\gamma
_{k,m}}-1}-\frac{3^{\gamma}-1}{2^{\gamma}-1} \biggr)
\\
&&\qquad=\frac{3^{\gamma}-1}{2^{\gamma}-1} \biggl[\sqrt{k} \biggl(\frac{3^{\hat
\gamma_{k,m}}-1}{3^{\gamma}-1}-1 \biggr)-
\sqrt{k} \biggl(\frac{2^{\hat\gamma_{k,m}}-1}{2^{\gamma}-1}-1 \biggr) \biggr]
\end{eqnarray*}
has the same limit distribution as
\[
\frac{3^{\gamma}-1}{2^{\gamma}-1}\sqrt{k} (\hat\gamma_{k,m}-\gamma ) \biggl(
\frac{\log3}{1-3^{-\gamma}}-\frac{\log2}{1-2^{-\gamma
}} \biggr)
\]
and, consequently,
\begin{eqnarray*}
&&\sqrt{k} (\hat\gamma_{k,m}-\gamma )
\\
&&\qquad\to^d \frac{1}{\Gamma(1-\gamma)} \biggl(\frac{\log3}{1-3^{-\gamma
}}-
\frac{\log2}{1-2^{-\gamma}} \biggr)^{-1} \\
&&\hspace*{8pt}\qquad\quad{}\times\biggl(\frac\gamma{3^{\gamma}-1}
(Q_2-Q_0 )-\frac\gamma {2^{\gamma}-1}
(Q_1-Q_0 ) \biggr). %&&\frac\lambda{\Gamma(1-\gamma)}\left(\frac{\log3}{1-3^{-\gamma}}-
\end{eqnarray*}

For the asymptotic distribution of $\hat a_{k,m}$ we write
\begin{eqnarray*}
&&\sqrt{k} \biggl(\frac{\hat a_{k,m}}{a_m}-1 \biggr)\frac{\hat\gamma
_{k,m}}{ (2^{\hat\gamma_{k,m}}-1 )\Gamma (1-\hat\gamma
_{k,m} )}
\\
&&\qquad{}\times\biggl\{\sqrt{k} \biggl(\frac{2\beta_1-\beta_0}{a_m}-\frac{2^\gamma
-1}\gamma\Gamma(1-
\gamma) \biggr)\\
&&\hspace*{8pt}\qquad\quad{}+ \sqrt{k} \biggl(\frac{2^\gamma-1}\gamma\Gamma(1-\gamma)-
\frac{2^{\hat
\gamma_{k,m}}-1}{\hat\gamma_{k,m}}\Gamma (1-\hat\gamma_{k,m} ) \biggr) \biggr\},
\end{eqnarray*}
and the statement follows, for example, by Cram\'er's delta method.

For the asymptotic distribution of $\hat b_{k,m}$, we write
\begin{eqnarray*}
&&\sqrt{k} \biggl(\frac{\hat b_{k,m}-b_m}{a_m} \biggr)\\
&&\qquad=\sqrt{k} \biggl(\frac
{\beta_0-b_m}{a_m}-
\frac{\Gamma(1-\gamma)-1}{\gamma} \biggr)
\\
&&\qquad\quad{}-\frac{\hat a_{k,m}}{a_m}\sqrt{k} \biggl(\frac{\Gamma(1-\hat\gamma
_{k,m})-1}{\hat\gamma_{k,m}}-\frac{\Gamma(1-\gamma)-1}{\gamma}
\biggr)\\
&&\qquad\quad{}+\frac{\Gamma(1-\gamma)-1}{\gamma}\sqrt{k} \biggl(\frac{\hat
a_{k,m}}{a_m}-1 \biggr)
\end{eqnarray*}
and the statement follows, for example, by Cram\'er's delta method.
\end{pf*}

\begin{pf*}{Proof of Theorem~\ref{AsymptNquantile}}
The proof follows the line of the comparable result for the POT method
[see, e.g., \citet{HaanFerreira}, Chapter~4.3]. Let
$c_n=1/(mp_n)$. Then
\begin{eqnarray*}
&&\frac{\sqrt k (\hat x_{k,m}-x_n )}{a_m q_\gamma
(c_n)}\\
&&\qquad= \frac{\sqrt k}{a_m q_\gamma(c_n)} \biggl(\hat b_{k,m}+\hat
a_{k,m}\frac{c_n^{\hat\gamma_{k,m}}-1}{\hat\gamma_{k,m}}-V \biggl(\frac{1}{-\log(1-p_n)} \biggr) \biggr)
\\
&&\qquad=\frac{\sqrt k}{q_\gamma(c_n)} \frac{\hat b_{k,m}-b_m}{a_m}+\frac
{\hat a_{k,m}}{a_m}
\frac{\sqrt k}{q_\gamma(c_n)} \biggl(\frac
{c_n^{\hat\gamma_{k,m}}-1}{\hat\gamma}-\frac{c_n^{\gamma}-1}{\gamma
} \biggr)
\\
&&\qquad\quad{}-\frac{\sqrt k}{q_\gamma(c_n)} \biggl(\frac{V ({m}/{(-m\log
(1-p_n))} )-V(m)}{a_m}-\frac{c_n^{\gamma}-1}{\gamma} \biggr)\\
&&\qquad\quad{}+
\frac
{c_n^{\gamma}-1}{\gamma q_\gamma(c_n)}\sqrt{k} \biggl(\frac{\hat
a_{k,m}}{a_m}-1 \biggr).
\end{eqnarray*}
Similarly, as on pages 135--137 of \citet{HaanFerreira}, this
converges in distribution to
\[
\Delta+(\gamma_-)^2 \Xi-\gamma_-\Lambda-\lambda
\frac{\gamma_-}{\gamma
_-+\rho}.
\]
\upqed\end{pf*}

\begin{appendix}
%s5 #&#
\section*{Appendix: Asymptotic variances and biases of the PWM estimators}

The following provides a basis for an algorithm to calculate the
asymptotic variances/covariances and biases of the PWM estimators.

Let $Q_r=(r+1)\int_0^1 s^{r-1}(-\log s)^{-1-\gamma}E(s) \,ds+\lambda
I_r(\gamma,\rho)$, $r=0,1,2$, as defined in Theorem~\ref{AsymptNMr}.
For $r,m=0,1,2$,
%
%e12 #&#
\begin{eqnarray}\label{covQr}
&&\operatorname{Cov}(Q_r,Q_m)\nonumber
\\
&&\qquad=(r+1) (m+1)\nonumber\\
&&\qquad\quad{}\times\int_0^1\int
_0^1 s^{r-1}u^{m-1}(-\log
s)^{-1-\gamma}(-\log u)^{-1-\gamma}EB(s)B(u) \,ds \,du
\nonumber
\\[-8pt]
\\[-8pt]
\nonumber
&&\qquad=(r+1) (m+1)\int_0^1 u^{m-1}(1-u)
(-\log u)^{-1-\gamma} \int_0^u
s^{r}(-\log s)^{-1-\gamma} \,ds \,du
\\
&& \qquad\quad{}+ (r+1) (m+1)\nonumber\\
&&\qquad\quad{}\times\int_0^1 s^{r-1}(1-s)
(-\log s)^{-1-\gamma}\int_0^s
u^{m}(-\log u)^{-1-\gamma} \,du \,ds\nonumber
\end{eqnarray}
using the fact that $EB(s)B(u)=\min(s,u)-su=s(1-u)$ for $0<s<u$.
These integrals can be evaluated numerically (we have used Mathematica
software). %Denote the elements of the correspondent
%variance/covariance matrix by $\sigma_{Q_0}^2$, $\sigma_{Q_{01}}$, $

From Theorems \ref{AsymptNgab} and \ref{AsymptNquantile}, after some
calculations,
%
%e13 #&#
%e14 #&#
%e15 #&#
%e16 #&#
\begin{eqnarray}
\sqrt k (\hat\gamma_{k,m}-\gamma )&\to^d&C_\gamma
(k_{\gamma,0}Q_0+k_{\gamma,1}Q_1+k_{\gamma,2}Q_2
),\label
{gQexpans}
\\
\sqrt k \biggl(\frac{\hat a_{k,m}}{a_m}-1 \biggr)&\to^d&
k_{a,0}Q_0+k_{a,1}Q_1+k_{a,2}Q_2,
\label{aQexpans}
\\
\sqrt k\frac{\hat b_{k,m}-b_m}{a_m}&\to^d&k_{b,0}Q_0+k_{b,1}Q_1+k_{b,2}Q_2,\label{bQexpans}
\\
\sqrt k\frac{ (\hat x_{k,m}-x_n )}{a_m q_\gamma(c_n)}&\to^d&k_{x,0}Q_0+k_{x,1}Q_1+k_{x,2}Q_2,\label{xQexpans}
\end{eqnarray}
where, for $\gamma\neq 0$,
\begin{eqnarray*}
C_\gamma&=&\frac{1}{\Gamma(1-\gamma)} \biggl(\frac{\log3}{1-3^{-\gamma
}}-
\frac{\log2}{1-2^{-\gamma}} \biggr)^{-1},
\\
k_{\gamma,0}&=&\frac{\gamma(3^{\gamma}-2^{\gamma})}{(3^{\gamma
}-1)(2^{\gamma}-1)},\qquad k_{\gamma,1}=\frac{-\gamma}{2^{\gamma
}-1},\qquad
k_{\gamma,2}=\frac{\gamma}{3^{\gamma}-1};
\\
C_a&=&\frac{\log2}{\gamma} \biggl(\frac{1}{\log2}-
\frac{\gamma
}{1-2^{-\gamma}} \biggr)+\frac{\Gamma'(1-\gamma)}{\Gamma(1-\gamma)},\\
 k_{a,0}&=&C_\gamma
k_{\gamma,0}C_a-\frac{\gamma}{(2^{\gamma}-1)\Gamma
(1-\gamma)},
\\
k_{a,1}&=&C_\gamma k_{\gamma,1}C_a+
\frac{\gamma}{(2^{\gamma}-1)\Gamma
(1-\gamma)},\\
 k_{a,2}&=&C_\gamma k_{\gamma,2}C_a-
\frac{\gamma
}{(2^{\gamma}-1)\Gamma(1-\gamma)};
\\
C_b&=&\frac{\gamma\Gamma'(1-\gamma)-1+\Gamma(1-\gamma)}{\gamma^2}, \\
 k_{b,0}&=&1+C_\gamma
k_{\gamma,0}C_b+k_{\gamma,0}\frac{1-\Gamma(1-\gamma
)}{\gamma},
\\
k_{b,1}&=&C_\gamma k_{\gamma,1}C_b+k_{\gamma,1}
\frac{1-\Gamma(1-\gamma
)}{\gamma},\\
 k_{b,2}&=&C_\gamma k_{\gamma,2}C_b+k_{\gamma,2}
\frac
{1-\Gamma(1-\gamma)}{\gamma};
\\
k_{x,0}&=&C_\gamma k_{\gamma,0}+(\gamma_-)^2k_{b,0}-
\gamma_-k_{a,0},\\
k_{x,1}&=&C_\gamma k_{\gamma,1}+(\gamma_-)^2k_{b,1}-
\gamma_-k_{a,1},\qquad k_{x,2}=C_\gamma k_{\gamma,2}+(
\gamma_-)^2k_{b,2}-\gamma_-k_{a,2}
\end{eqnarray*}
and, for $\gamma=0$,
\begin{eqnarray*}
C_\gamma&=&2\bigl(\log(3/2)\bigr)^{-1},\qquad k_{\gamma,0}=(
\log2)^{-1}-(\log 3)^{-1}, \\
 k_{\gamma,1}&=&-(
\log2)^{-1},\qquad k_{\gamma,2}=(\log3)^{-1};
%
%k_{a,0}=\frac{\log(2/3)+2\Gamma'(1)}{\log2\log3}  k_{a,1}=\frac{
%k_{a,2}=\frac{2\log3}{\log3-\log2}\left(\frac{\log2}2+\Gamma'(1)
\\
C_a&=&2^{-1}\log2+\Gamma'(1),\qquad
C_b=-\Gamma''(1)
\end{eqnarray*}
and the rest follow similarly by continuity. Then the asymptotic
variances, covariances and biases follow by combining \eqref{covQr}
with \eqref{gQexpans}--\eqref{xQexpans} in the obvious way.
%Explicit (complicated) expressions for the limiting covariance matrix
%can be found in \citet{HoskingWallisWood85}, cf. $v_{r,r}$,
%$v_{r,r+1}$ and $v_{r,r+s}$ (C.9)--(C.11) in their Appendix C. From
%there, var$(Q_r)=(r+1)^2 v_{rr}$ and cov$(Q_r,Q_s)=(r+1)(s+1)v_{r,s}$.
\end{appendix}

\section*{Acknowledgements}
We thank Holger Drees for a useful suggestion.

We would like to thank three unknown referees for their genuine
interest and their insightful comments.

% imsref loaded by akundreckaite, 2014-11-06 08:10:52


%suskaldyti doi

\printaddresses

\begin{thebibliography}{26}
% pybtex-1.20. Style name=ims, version=2.91, label_style=nameyear, sorting_style=complex, cfg=None, language=None.


%b1 ###
%b1 #&#
\bibitem[\protect\citeauthoryear{B{\"u}cher and Segers}{2014}]{BucherSegers14}
\begin{barticle}[mr]
\bauthor{\bsnm{B{\"u}cher},~\bfnm{Axel}\binits{A.}} \AND
\bauthor{\bsnm{Segers},~\bfnm{Johan}\binits{J.}}
(\byear{2014}).
\btitle{Extreme value copula estimation based on block maxima of a multivariate stationary time series}.
\bjournal{Extremes}
\bvolume{17}
\bpages{495--528}.
\bid{doi={10.1007/s10687-014-0195-8}, issn={1386-1999}, mr={3252823}}
\end{barticle}
%
%\OrigBibText
%B\"ucher, A. and Segers, J. (2014) Extreme
%value copula estimation based on block maxima of a multivariate
%stationary time series. To appear in {\it Extremes}, arXiv:1311.3060.
%\endOrigBibText
\bptok{imsref}%
% NOT OUTPUTTED:
%   number = 3
%   doi = http://dx.doi.org/10.1007/s10687-014-0195-8
%   fjournal = Extremes. Statistical Theory and Applications in Science, Engineering and Economics
\endbibitem

%b2 ###
%b2 #&#
\bibitem[\protect\citeauthoryear{Cai, de~Haan and Zhou}{2013}]{CaiHaanZhou}
\begin{barticle}[mr]
\bauthor{\bsnm{Cai},~\bfnm{Juan-Juan}\binits{J.-J.}},
\bauthor{\bparticle{de} \bsnm{Haan},~\bfnm{Laurens}\binits{L.}} \AND
\bauthor{\bsnm{Zhou},~\bfnm{Chen}\binits{C.}}
(\byear{2013}).
\btitle{Bias correction in extreme value statistics with index around zero}.
\bjournal{Extremes}
\bvolume{16}
\bpages{173--201}.
\bid{doi={10.1007/s10687-012-0158-x}, issn={1386-1999}, mr={3057195}}
\end{barticle}
%
%\OrigBibText
%Cai, J.-J., de Haan, L. and Zhou, C. (2013) Bias
%correction in extreme value statistics with index around zero. {\it
%Extremes} {\bf16}, 173--201.
%\endOrigBibText
\bptok{imsref}%
% NOT OUTPUTTED:
%   number = 2
%   doi = http://dx.doi.org/10.1007/s10687-012-0158-x
%   fjournal = Extremes. Statistical Theory and Applications in Science, Engineering and Economics
\endbibitem

%b3 ###
%b3 #&#
\bibitem[\protect\citeauthoryear{Caires}{2009}]{Caires09}
\begin{bmisc}[auto:parserefs-M02]
\bauthor{\bsnm{Caires},~\bfnm{S.}\binits{S.}}
(\byear{2009}).
\bhowpublished{A comparative simulation study of the annual maxima
and the peaks-over-threshold methods. SBW-Belastingen: Subproject ``Statistics''.
Deltares Report 1200264-002}.
\end{bmisc}
%
%\OrigBibText
%Caires, S. (2009) A comparative simulation study of the annual maxima
%and the peaks-over-threshold methods. SBW-Belastingen: subproject
%`Statistics'. Deltares Report 1200264-002.
%\endOrigBibText
\bptok{imsref}%
\endbibitem

%b4 ###
%b4 #&#
\bibitem[\protect\citeauthoryear{Cs{\"o}rg{\H{o}} and Horv{\'a}th}{1993}]{Csorgo}
\begin{bbook}[mr]
\bauthor{\bsnm{Cs{\"o}rg{\H{o}}},~\bfnm{Mikl{\'o}s}\binits{M.}} \AND
\bauthor{\bsnm{Horv{\'a}th},~\bfnm{Lajos}\binits{L.}}
(\byear{1993}).
\btitle{Weighted Approximations in Probability and Statistics}.
%Probability and Mathematical Statistics}.
\bpublisher{Wiley},
\blocation{Chichester}.
\bid{mr={1215046}}
\end{bbook}
%
%\OrigBibText
%Cs\"org\H o, M. and Horv\'ath, L. (1993) Weighted Approximations in
%Probability and Statistics.
%John Wiley \& Sons, Chichester, England
%\endOrigBibText
\bptok{imsref}%
% NOT OUTPUTTED:
%   isbn = 0-471-93635-9
%   fpage = xvi+442
\endbibitem

%b5 ###
%b5 #&#
\bibitem[\protect\citeauthoryear{Cunnane}{1973}]{Cunnane73}
\begin{barticle}[auto:parserefs-M02]
\bauthor{\bsnm{Cunnane},~\bfnm{C.}\binits{C.}}
(\byear{1973}).
\btitle{A particular comparison of annual maxima and partial
duration series methods of flood frequency prediction}.
\bjournal{J. Hydrol.}
\bvolume{18}
\bpages{257--271}.
\end{barticle}
%
%\OrigBibText
%Cunnane, C. (1973) A particular comparison of
%annual maxima and partial duration series methods of flood frequency
%prediction. {\it J. Hydrol.} {\bf18}, 257--271.
%\endOrigBibText
\bptok{imsref}%
\endbibitem

%b6 ###
%b6 #&#
\bibitem[\protect\citeauthoryear{de~Haan and Ferreira}{2006}]{HaanFerreira}
\begin{bbook}[mr]
\bauthor{\bparticle{de} \bsnm{Haan},~\bfnm{Laurens}\binits{L.}} \AND
\bauthor{\bsnm{Ferreira},~\bfnm{Ana}\binits{A.}}
(\byear{2006}).
\btitle{Extreme Value Theory: An Introduction}.
\bpublisher{Springer},
\blocation{New York}.
\bid{doi={10.1007/0-387-34471-3}, mr={2234156}}
\end{bbook}
%
%\OrigBibText
%de Haan, L. and Ferreira, A. (2006) Extreme
%Value Theory: An Introduction. Springer, Boston.
%\endOrigBibText
\bptok{imsref}%
% NOT OUTPUTTED:
%   doi = http://dx.doi.org/10.1007/0-387-34471-3
%   isbn = 978-0-387-23946-0; 0-387-23946-4
%   fpage = xviii+417
\endbibitem

%b7 ###
%b7 #&#
\bibitem[\protect\citeauthoryear{de~Valk}{1993}]{Valk93}
\begin{bmisc}[auto:parserefs-M02]
\bauthor{\bparticle{de} \bsnm{Valk},~\bfnm{C.}\binits{C.}}
(\byear{1993}).
\bhowpublished{Estimation of marginals from measurements and hindcast data.
WL{|}Delft Hydraulics Report H1700}.
\end{bmisc}
%
%\OrigBibText
%de Valk, C. (1993) Estimation of marginals from measurements and
%hindcast data. WL|Delft Hydraulics Report H1700.
%\endOrigBibText
\bptok{imsref}%
\endbibitem

%b8 ###
%b8 #&#
\bibitem[\protect\citeauthoryear{Diebolt et~al.}{2008}]{DieboltGuillouNaveauRibereau08}
\begin{barticle}[mr]
\bauthor{\bsnm{Diebolt},~\bfnm{Jean}\binits{J.}},
\bauthor{\bsnm{Guillou},~\bfnm{Armelle}\binits{A.}},
\bauthor{\bsnm{Naveau},~\bfnm{Philippe}\binits{P.}} \AND
\bauthor{\bsnm{Ribereau},~\bfnm{Pierre}\binits{P.}}
(\byear{2008}).
\btitle{Improving probability-weighted moment methods for the generalized extreme value distribution}.
\bjournal{REVSTAT}
\bvolume{6}
\bpages{35--50}.
\bid{issn={1645-6726}, mr={2386298}}
\end{barticle}
%
%\OrigBibText
%Diebolt, J., Guillou, A., Naveau, P. and Ribereau, P. (2008) Improving
%probability-weighted moment methods for the generalized extreme value
%distribution. {\it RevStat} {\bf6}, 33--50.
%\endOrigBibText
\bptok{imsref}%
% NOT OUTPUTTED:
%   number = 1
%   fjournal = REVSTAT Statistical Journal
\endbibitem

%b9 ###
%b9 #&#
\bibitem[\protect\citeauthoryear{Dombry}{2013}]{Dombry}
\begin{bmisc}[auto:parserefs-M02]
\bauthor{\bsnm{Dombry},~\bfnm{C.}\binits{C.}}
(\byear{2013}).
\bhowpublished{Maximum likelihood estimators for the extreme value
index based on the block maxima method.
Available at \arxivurl{arXiv:1301.5611}}.
\end{bmisc}
%
%\OrigBibText
%Dombry, C. (2013) Maximum likelihood estimators for
%the extreme value index based on the block maxima method: arXiv:1301.5611
%\endOrigBibText
\bptok{imsref}%
\endbibitem\vadjust{\goodbreak}

%b10 ###
%b10 #&#
\bibitem[\protect\citeauthoryear{Drees}{1998}]{Drees98}
\begin{barticle}[mr]
\bauthor{\bsnm{Drees},~\bfnm{Holger}\binits{H.}}
(\byear{1998}).
\btitle{On smooth statistical tail functionals}.
\bjournal{Scand. J. Stat.}
\bvolume{25}
\bpages{187--210}.
\bid{doi={10.1111/1467-9469.00097}, issn={0303-6898}, mr={1614276}}
\end{barticle}
%
%\OrigBibText
%Drees, H. (1998) On smooth statistical tail
%functionals. {\it Scand. J. Statist.} \textbf{25}, 187--210.
%\endOrigBibText
\bptok{imsref}%
% NOT OUTPUTTED:
%   number = 1
%   doi = http://dx.doi.org/10.1111/1467-9469.00097
%   fjournal = Scandinavian Journal of Statistics. Theory and Applications
\endbibitem

%b11 ###
%b11 #&#
\bibitem[\protect\citeauthoryear{Drees, de~Haan and Li}{2003}]{DreesHaanLi03}
\begin{barticle}[mr]
\bauthor{\bsnm{Drees},~\bfnm{Holger}\binits{H.}},
\bauthor{\bparticle{de} \bsnm{Haan},~\bfnm{Laurens}\binits{L.}} \AND
\bauthor{\bsnm{Li},~\bfnm{Deyuan}\binits{D.}}
(\byear{2003}).
\btitle{On large deviation for extremes}.
\bjournal{Statist. Probab. Lett.}
\bvolume{64}
\bpages{51--62}.
\bid{doi={10.1016/S0167-7152(03)00130-5}, issn={0167-7152}, mr={1995809}}
\end{barticle}
%
%\OrigBibText
%Drees, H., de Haan, L. and Li, D. (2003) On
%large deviation for extremes. {\it Statistics \& Probability Letters}
%\textbf{64}, 51--62.
%\endOrigBibText
\bptok{imsref}%
% NOT OUTPUTTED:
%   number = 1
%   doi = http://dx.doi.org/10.1016/S0167-7152(03)00130-5
%   coden = SPLTDC
%   fjournal = Statistics \& Probability Letters
\endbibitem

%b12 ###
%b12 #&#
\bibitem[\protect\citeauthoryear{Gumbel}{1958}]{Gumbel58}
\begin{bbook}[mr]
\bauthor{\bsnm{Gumbel},~\bfnm{E.~J.}\binits{E.~J.}}
(\byear{1958}).
\btitle{Statistics of Extremes}.
\bpublisher{Columbia Univ. Press},
\blocation{New York}.
\bid{mr={0096342}}
\end{bbook}
%
%\OrigBibText
%Gumbel, E. (1958) Statistics of Extremes, Columbia
%University Press.
%\endOrigBibText
\bptok{imsref}%
% NOT OUTPUTTED:
%   fpage = xx+375
\endbibitem

%b13 ###
%b13 #&#
\bibitem[\protect\citeauthoryear{Hosking}{1990}]{Hosking90}
\begin{barticle}[mr]
\bauthor{\bsnm{Hosking},~\bfnm{J.~R.~M.}\binits{J.~R.~M.}}
(\byear{1990}).
\btitle{{$L$}-moments: Analysis and estimation of distributions using linear combinations of order statistics}.
\bjournal{J. R. Stat. Soc. Ser. B Stat. Methodol.}
\bvolume{52}
\bpages{105--124}.
\bid{issn={0035-9246}, mr={1049304}}
\end{barticle}
%
%\OrigBibText
%Hosking, J.R.M. (1990) L-Moments: Analysis and
%estimation of distributions using linear combinations of order
%statistics. {J. R. Stat. Soc. B} \textbf{52}, 105--124.
%\endOrigBibText
\bptok{imsref}%
% NOT OUTPUTTED:
%   url = http://links.jstor.org/sici?sici=0035-9246(1990)52:1<105:AAEODU>2.0.CO;2-7&origin=MSN
%   number = 1
%   coden = JSTBAJ
%   fjournal = Journal of the Royal Statistical Society. Series B. Methodological
\endbibitem

%b14 ###
%b14 #&#
\bibitem[\protect\citeauthoryear{Hosking and Wallis}{1987}]{HoskingWallis87}
\begin{barticle}[mr]
\bauthor{\bsnm{Hosking},~\bfnm{J.~R.~M.}\binits{J.~R.~M.}} \AND
\bauthor{\bsnm{Wallis},~\bfnm{J.~R.}\binits{J.~R.}}
(\byear{1987}).
\btitle{Parameter and quantile estimation for the generalized {P}areto distribution}.
\bjournal{Technometrics}
\bvolume{29}
\bpages{339--349}.
\bid{doi={10.2307/1269343}, issn={0040-1706}, mr={0906643}}
\end{barticle}
%
%\OrigBibText
%Hosking, J.R.M. and Wallis, J.R. (1987)
%Parameter and quantile estimation for the Generalized Pareto
%Distribution. {\it Technometrics} \textbf{29}, 339--349.
%\endOrigBibText
\bptok{imsref}%
% NOT OUTPUTTED:
%   number = 3
%   doi = http://dx.doi.org/10.2307/1269343
%   coden = TCMTA2
%   fjournal = Technometrics. A Journal of Statistics for the Physical, Chemical and Engineering Sciences
\endbibitem

%b15 ###
%b15 #&#
\bibitem[\protect\citeauthoryear{Hosking, Wallis and Wood}{1985}]{HoskingWallisWood85}
\begin{barticle}[mr]
\bauthor{\bsnm{Hosking},~\bfnm{J.~R.~M.}\binits{J.~R.~M.}},
\bauthor{\bsnm{Wallis},~\bfnm{J.~R.}\binits{J.~R.}} \AND
\bauthor{\bsnm{Wood},~\bfnm{E.~F.}\binits{E.~F.}}
(\byear{1985}).
\btitle{Estimation of the generalized extreme-value distribution by the method of probability-weighted moments}.
\bjournal{Technometrics}
\bvolume{27}
\bpages{251--261}.
\bid{doi={10.2307/1269706}, issn={0040-1706}, mr={0797563}}
\end{barticle}
%
%\OrigBibText
%Hosking, J.R.M., Wallis, J.R. and Wood,
%E.F. (1985) Estimation of the Generalized Extreme-Value Distribution by
%the Method of Probability Weighted Moments. {\it Technometrics} \textbf
%{27}, 251--261.
%\endOrigBibText
\bptok{imsref}%
% NOT OUTPUTTED:
%   number = 3
%   doi = http://dx.doi.org/10.2307/1269706
%   coden = TCMTA2
%   fjournal = Technometrics. A Journal of Statistics for the Physical, Chemical and Engineering Sciences
\endbibitem

%b16 ###
%b16 #&#
\bibitem[\protect\citeauthoryear{Katz, Parlange and Naveau}{2002}]{KatzParlangeNaveau02}
\begin{barticle}[auto:parserefs-M02]
\bauthor{\bsnm{Katz},~\bfnm{R.~W.}\binits{R.~W.}}
\bauthor{\bsnm{Parlange},~\bfnm{M.~B.}\binits{M.~B.}}
\AND
\bauthor{\bsnm{Naveau},~\bfnm{P.}\binits{P.}}
(\byear{2002}).
\btitle{Statistics of extremes in hydrology}.
\bjournal{Advances in Water Resources}
\bvolume{25}
\bpages{1287--1304}.
\end{barticle}
%
%\OrigBibText
%RW Katz, MB , P
%Katz, Parlange and Naveau (2002) Statistics of extremes in hydrology
%{\it Advances in Water Resources} {\bf25}, 1287--1304.
%\endOrigBibText
\bptok{imsref}%
\endbibitem

%b17 ###
%b17 #&#
\bibitem[\protect\citeauthoryear{Kharin et~al.}{2007}]{Kharin07}
\begin{barticle}[auto:parserefs-M02]
\bauthor{\bsnm{Kharin},~\bfnm{V.~V.}\binits{V.~V.}},
\bauthor{\bsnm{Zwiers},~\bfnm{F.~W.}\binits{F.~W.}},
\bauthor{\bsnm{Zhang},~\bfnm{X.}\binits{X.}} \AND
\bauthor{\bsnm{Hegerl},~\bfnm{G.~C.}\binits{G.~C.}}
(\byear{2007}).
\btitle{Changes in temperature and precipitation extremes in the IPCC
ensemble of global coupled model simulations}.
\bjournal{Journal of Climate}
\bvolume{20}
\bpages{1419--1444}.
\end{barticle}
%
%\OrigBibText
%Kharin, V.V., Zwiers, F.W., Zhang, X. and Hegerl,
%G.C. (2007) Changes in temperature and precipitation extremes in the
%IPCC ensemble of gloubal coupled model simulations. {\it Journal of
%Climate} {\bf20}, 1419--1444.
%\endOrigBibText
\bptok{imsref}%
\endbibitem

%b18 ###
%b18 #&#
\bibitem[\protect\citeauthoryear{Landwehr, Matalas and Wallis}{1979}]{LandwehrMatalasWallis79}
\begin{barticle}[auto:parserefs-M02]
\bauthor{\bsnm{Landwehr},~\bfnm{J.}\binits{J.}},
\bauthor{\bsnm{Matalas},~\bfnm{N.}\binits{N.}} \AND
\bauthor{\bsnm{Wallis},~\bfnm{J.}\binits{J.}}
(\byear{1979}).
\btitle{Probability weighted moments compared with some traditional
techniques in estimating Gumbel parameters and quantiles}.
\bjournal{Water Resources Research}
\bvolume{15}
\bpages{1055--1064}.
\end{barticle}
%
%\OrigBibText
%Landwehr, J., Matalas, N. and Wallis,
%J. (1979) Probability weighted moments compared with some traditional
%techniques in estimating Gumbel parameters and quantiles. {\it Water
%Resources Research} {\bf15}, 1055--1064.
%\endOrigBibText
\bptok{imsref}%
\endbibitem

%b19 ###
%b19 #&#
\bibitem[\protect\citeauthoryear{Madsen, Pearson and Rosbjerg}{1997}]{MadsenPearsonRosbjerg97}
\begin{barticle}[auto:parserefs-M02]
\bauthor{\bsnm{Madsen},~\bfnm{H.}\binits{H.}},
\bauthor{\bsnm{Pearson},~\bfnm{C.~P.}\binits{C.~P.}} \AND
\bauthor{\bsnm{Rosbjerg},~\bfnm{D.}\binits{D.}}
(\byear{1997}).
\btitle{Comparison of annual maximum series and partial duration
series methods for modeling extreme hydrologic events 2. Regional modeling}.
\bjournal{Water Resources Research}
\bvolume{33}
\bpages{759--769}.
\end{barticle}
%
%\OrigBibText
%Madsen, H., Pearson, C.P. and
%Rosbjerg, D. (1997)
%Comparison of annual maximum series and partial duration series methods
%for modeling extreme hydrologic events 2. Regional modeling.
%{\it Water Resources Research} {\bf33}, 759--769.
%\endOrigBibText
\bptok{imsref}%
\endbibitem

%b20 ###
%b20 #&#
\bibitem[\protect\citeauthoryear{Madsen, Rasmussen and Rosbjerg}{1997}]{MadsenRasmussenRosbjerg97}
\begin{barticle}[auto:parserefs-M02]
\bauthor{\bsnm{Madsen},~\bfnm{H.}\binits{H.}},
\bauthor{\bsnm{Rasmussen},~\bfnm{P.~F.}\binits{P.~F.}} \AND
\bauthor{\bsnm{Rosbjerg},~\bfnm{D.}\binits{D.}}
(\byear{1997}).
\btitle{Comparison of annual maximum series and partial duration
series methods for modeling extreme hydrologic events 1. At-site modeling}.
\bjournal{Water Resources Research}
\bvolume{33}
\bpages{747--757}.
\end{barticle}
%
%\OrigBibText
%Madsen, H., Rasmussen, P.F. and Rosbjerg, D. (1997) Comparison of
%annual maximum series and partial duration series methods for modeling
%extreme hydrologic events 1. At-site modeling. {\it Water Resources
%Research} {\bf33}, 747--757.
%\endOrigBibText
\bptok{imsref}%
\endbibitem

%b21 ###
%b21 #&#
\bibitem[\protect\citeauthoryear{Martins and Stedinger}{2001}]{MartinsStedinger}
\begin{barticle}[auto:parserefs-M02]
\bauthor{\bsnm{Martins},~\bfnm{E.~S.}\binits{E.~S.}} \AND
\bauthor{\bsnm{Stedinger},~\bfnm{J.~R.}\binits{J.~R.}}
(\byear{2001}).
\btitle{Historical information in a generalized maximum likelihood
framework with partial duration and annual maximum series}.
\bjournal{Water Resources Research}
\bvolume{37}
\bpages{2559--2567}.
\end{barticle}
%
%\OrigBibText
%Martins, E.S. and Stedinger, J.R. (2001)
%Historical information in a generalized maximum likelihood framework
%with partial duration and annual maximum series. {\it Water Resources
%Research} {\bf37}, 2559--2567.
%\endOrigBibText
\bptok{imsref}%
\endbibitem

%b22 ###
%b22 #&#
\bibitem[\protect\citeauthoryear{Naveau et~al.}{2009}]{NaveauGuillouCooleyDiebolt09}
\begin{barticle}[mr]
\bauthor{\bsnm{Naveau},~\bfnm{Philippe}\binits{P.}},
\bauthor{\bsnm{Guillou},~\bfnm{Armelle}\binits{A.}},
\bauthor{\bsnm{Cooley},~\bfnm{Daniel}\binits{D.}} \AND
\bauthor{\bsnm{Diebolt},~\bfnm{Jean}\binits{J.}}
(\byear{2009}).
\btitle{Modelling pairwise dependence of maxima in space}.
\bjournal{Biometrika}
\bvolume{96}
\bpages{1--17}.
\bid{doi={10.1093/biomet/asp001}, issn={0006-3444}, mr={2482131}}
\end{barticle}
%
%\OrigBibText
%Naveau, P., Guillou, A., Cooley, D. and Diebolt, J. (2009) Modelling
%pairwise dependence of maxima in space. {\it Biometrika} {\bf96}, 1--17.
%\endOrigBibText
\bptok{imsref}%
% NOT OUTPUTTED:
%   number = 1
%   doi = http://dx.doi.org/10.1093/biomet/asp001
%   coden = BIOKAX
%   fjournal = Biometrika
\endbibitem

%b23 ###
%b23 #&#
\bibitem[\protect\citeauthoryear{Pickands}{1975}]{Pickands75}
\begin{barticle}[mr]
\bauthor{\bsnm{Pickands},~\bfnm{James}\binits{J.} \bsuffix{III}}
(\byear{1975}).
\btitle{Statistical inference using extreme order statistics}.
\bjournal{Ann. Statist.}
\bvolume{3}
\bpages{119--131}.
\bid{issn={0090-5364}, mr={0423667}}
\end{barticle}
%
%\OrigBibText
%Pickands, J. III (1975) Statistical inference using extreme order
%statistics. {\it Ann. Statist.} {\bf3}, 119--131.
%\endOrigBibText
\bptok{imsref}%
% NOT OUTPUTTED:
%   fjournal = The Annals of Statistics
\endbibitem

%b24 ###
%b24 #&#
\bibitem[\protect\citeauthoryear{Prescott and Walden}{1980}]{PrescottWalden1980}
\begin{barticle}[mr]
\bauthor{\bsnm{Prescott},~\bfnm{P.}\binits{P.}} \AND
\bauthor{\bsnm{Walden},~\bfnm{A.~T.}\binits{A.~T.}}
(\byear{1980}).
\btitle{Maximum likelihood estimation of the parameters of the generalized extreme-value distribution}.
\bjournal{Biometrika}
\bvolume{67}
\bpages{723--724}.
\bid{doi={10.1093/biomet/67.3.723}, issn={0006-3444}, mr={0601119}}
\end{barticle}
%
%\OrigBibText
%Prescott, P. and Walden, A.T. (1980) Maximum likelihood estimation of
%the parameters of the generalized extreme-value distribution. {\it
%Biometrika} {\bf67}, 723--724.
%\endOrigBibText
\bptok{imsref}%
% NOT OUTPUTTED:
%   number = 3
%   doi = http://dx.doi.org/10.1093/biomet/67.3.723
%   coden = BIOKAX
%   fjournal = Biometrika
\endbibitem

%b25 ###
%b25 #&#
\bibitem[\protect\citeauthoryear{van~den Brink, K{\"{o}}nnen and Opsteegh}{2005}]{BrinkKonnenOpsteegh05}
\begin{barticle}[pbm]
\bauthor{\bsnm{van~den Brink},~\bfnm{H.~W.}\binits{H.~W.}},
\bauthor{\bsnm{K{\"{o}}nnen},~\bfnm{G.~P.}\binits{G.~P.}} \AND
\bauthor{\bsnm{Opsteegh},~\bfnm{J.~D.}\binits{J.~D.}}
(\byear{2005}).
\btitle{Uncertainties in extreme surge level estimates from observational records}.
\bjournal{Philos. Trans. R. Soc. Lond. Ser. A Math. Phys. Eng. Sci.}
\bvolume{363}
\bpages{1377--1386}.
\bid{doi={10.1098/rsta.2005.1573}, issn={1364-503X}, pii={LL2HPE9025M4UAPF}, pmid={16191655}}
\end{barticle}
%
%\OrigBibText
%van den Brink, H.W., Konnen, G.P. and Opsteegh, J.D. (2005)
%Uncertainties in extreme surge level estimates from observational
%records. {\it Phil. Trans. R. Soc. A} {\bf363}, 1377--1386.
%\endOrigBibText
\bptok{imsref}%
% NOT OUTPUTTED:
%   number = 1831
%   fjournal = Philosophical transactions. Series A, Mathematical, physical, and engineering sciences
\endbibitem

%b26 ###
%b26 #&#
\bibitem[\protect\citeauthoryear{Wang}{1991}]{Wang91}
\begin{barticle}[auto:parserefs-M02]
\bauthor{\bsnm{Wang},~\bfnm{Q.~J.}\binits{Q.~J.}}
(\byear{1991}).
\btitle{The POT model described by the generalized Pareto distribution
with Poisson arrival rate}.
\bjournal{Journal of Hidrology}
\bvolume{129}
\bpages{263--280}.
\end{barticle}
%
%\OrigBibText
%Wang, Q.J. (1991) The POT model described by the generalized Pareto
%distribution with Poisson arrival rate. {\it Journal of Hidrology} {\bf
%129}, 263--280.
%\endOrigBibText
\bptok{imsref}%
\endbibitem
\end{thebibliography}
\end{document}